\documentclass[moor]{informs2}                   

\newcommand{\Rb}{\mathbb{R}}  \newcommand{\Nb}{\mathbb{N}}
\newcommand{\CS}{\mathcal{S}}
\newcommand{\upsub}{\mathop{\overline{\ominus}}}
\newcommand{\losub}{\mathop{\underline{\ominus}}}
\newcommand{\vep}{\varepsilon}
\newcommand{\DOP}{\mbox{OP}}

\usepackage{natbib}
 \NatBibNumeric
 \def\bibfont{\small}%
 \bibpunct[, ]{[}{]}{,}{n}{}{,}%

\TheoremsNumberedThrough     

\EquationsNumberedBySection 

\begin{document}


\RUNAUTHOR{Chistyakov and Pardalos}

\RUNTITLE{Concepts of Stability in Discrete Optimization
  Involving Generalized Addition Operations}

\TITLE{Concepts of Stability in Discrete Optimization\\
  Involving Generalized Addition Operations}

\ARTICLEAUTHORS{%
\AUTHOR{Vyacheslav V.~Chistyakov}
\AFF{Department of Applied Mathematics and Computer Science,
National Research University Higher School of Economics,
Bol'shaya Pech{\"e}rskaya Street 25/12, Nizhny Novgorod 603155,
Russian Federation, \EMAIL{vchistyakov@hse.ru}} 
\AUTHOR{Panos M.~Pardalos}
\AFF{Department of Industrial and Systems Engineering,
Center for Applied Optimization, University of Florida, 303 Weil Hall,
P.O.Box 116595, Gainesville, FL 32611-6595, USA, \EMAIL{pardalos@ufl.edu}} 
} 

\ABSTRACT{%
The paper addresses the tolerance approach to the sensitivity analysis of optimal
solutions to the nonlinear optimization problem of the form
  $$\mbox{$\bigoplus\limits_{y\in S}C(y)\to\min$\quad over\quad $S\in\CS$,}$$
 where $\CS$ is a collection of nonempty subsets of a finite set $X$
such that $\cup\CS=X$ and $\cap\CS=\varnothing$, $C$ is a cost (or weight)
function from $X$ into $\Rb^+=[0,\infty)$ or $(0,\infty)$, and $\oplus$ is a
continuous, associative, commutative, nondecreasing and unbounded binary operation
of generalized addition on $\Rb^+$, called an {\sf A}-operation. We evaluate and
present sharp estimates for upper and lower bounds of costs of elements from $X$,
for which an optimal solution to the above problem remains stable. These bounds
present new results in the sensitivity analysis as well as extend most known results in a unified way. We define an invariant of the optimization problem---the tolerance function,
which is independent of optimal solutions, and establish its basic properties, among which
we mention a characterization of the set of all optimal solutions, the uniqueness of optimal
solutions and extremal values of the tolerance function on an optimal solution.
}%


\MSCCLASS{Primary: 90C31; secondary: 90C27, 90C26}

\maketitle

%

\section{Introduction.} \label{s:intro}

The purpose of this paper is to introduce and study certain concepts of stability of
optimal solutions to the following nonlinear problem of discrete optimization:
  \begin{equation} \label{e:DO}
f_C(S)\equiv f(C)(S)=\bigoplus_{y\in S}C(y)\to\min,\qquad S\in\CS,
  \end{equation}
where $\CS$ is a collection of nonempty subsets (called trajectories) of a finite set
$X$ of cardinality $|X|\ge2$ such that $\cup\CS=X$ and $\cap\CS=\varnothing$
and $C:X\to\Rb^+$ is a given cost (or weight) function of elements from~$X$
with $\Rb^+=[0,\infty)$ or $(0,\infty)$. The objective function $f_C:\CS\to\Rb^+$
in \eqref{e:DO} is given by means of an operation $\oplus$ on the set $\Rb^+$, called
an {\sf A}-operation, which generalizes simultaneously the addition operation on
$\Rb^+=[0,\infty)$ and the operation of multiplication on $\Rb^+=(0,\infty)$.
More specifically, we assume that the operation $(u,v)\mapsto u\oplus v$ from
$\Rb^+\times\Rb^+$ into $\Rb^+$ is associative, commutative, nondecreasing
in each variable, unbounded in the sense that $u\oplus v\to\infty$ as $u\to\infty$
for all $v\in\Rb^+$, and continuous as a function of two real variables.
Simple examples of {\sf A}-operations are the usual addition
operation $u\oplus v=u+v$, the operation of taking the maximum
$u\oplus v=\max\{u,v\}$ and the usual operation of multiplication $u\oplus v=u\cdot v$.

In this paper we adopt the tolerance approach to the sensitivity analysis of optimal
solutions to problem \eqref{e:DO}: given an optimal solution $S^*\in\CS$ to problem
\eqref{e:DO} and an element $x\in X$, we are interested to what extent the cost $C(x)$
can be changed (the other costs remaining unperturbed) so that $S^*$ is still an
optimal solution. In other words, this can be expressed as how tolerant is the optimal
solution $S^*$ with respect to a change of the single cost $C(x)$. The tolerance
approach has been applied in the literature for a variety of combinatorial optimization
problems: shortest path and network flow problems \cite{Gus,Shier,Tarjan}, the
traveling salesman problem \cite{Jager,Libura,Sotskov,Turken}, the matrix coefficients
in linear programming problems \cite{Ravi}, the bottleneck problems
\cite{Chakra,Gal,GJM,Gord,Gordeev,Rama,Sotskov}. We refer to \cite{Gal,Greenberg,van}
for a comprehensive literature on the sensitivity (or post-solution) analysis.

Up to now almost all works in the sensitivity analysis, including the ones referred to above,
are concerned with two (most popular and natural) operations $\oplus$ in~\eqref{e:DO},
namely, the addition $+$ or $\max$, and so, the corresponding objective function
$f_C$ is either linear $f_C(S)=\sum_{y\in S}C(y)$ or bottleneck
$f_C(S)=\max_{y\in S}C(y)$, respectively. One may intuitively feel the difference
between the two operations: all costs contribute to the linear objective function equally
well, while only the largest costs contribute to the bottleneck objective function.
A more relevant explanation (in accordance with the theory of {\sf A}-operations
developed in Section~\ref{s:aop}) is that, given $v\ge0$, the function $\Phi(u)=u+v$
is strictly increasing (and so, the operation $+$ is termed to be strict), while the function
$\Psi(u)=\max\{u,v\}$ is only nondecreasing on $[0,\infty)$ (and so, the operation
$\max$ is not strict). In calculating the upper and lower bounds for costs of elements
from $X$, for which an optimal solution to problem \eqref{e:DO} remains stable,
one encounters the (generally, tacit) necessity to take the inverse(s) of the operation
under consideration (in our case $+$ or $\max$, cf. references above), which results in
inverting functions $\Phi$ and $\Psi$ as above. Since $\Phi$ is strictly increasing
(and continuous), its inverse, the subtraction, behaves well and causes no problem.
However, the usual inverse of the nondecreasing function $\Psi$ does not exist,
and so, one has to deal with its right and/or left inverse \cite{KR}, which is
why optimization problems with bottleneck objective functions are more complicated.
Technically speaking, problems \eqref{e:DO} with nonstrict {\sf A}-operations $\oplus$
tend to exhibit certain pathological, hysteretic-like properties involving retardness
(cf.~\cite{KPo}).
 
Given a generic {\sf A}-operation $\oplus$ on $\Rb^+$, we introduce two inverses
of~$\oplus$, called the upper subtraction $\upsub$ and the lower subtraction $\losub$.
It turns out that the two subtractions replace the ordinary subtraction (in the case of~$+$),
and they coincide on the common part of their domains iff ($=$\,if and only if) the
{\sf A}-operation $\oplus$ is strict. By means of the upper (lower) subtraction we
determine the upper (lower, respectively) stability interval of costs for an optimal
solution to problem~\eqref{e:DO}. In the case of strict {\sf A}-operations $\oplus$
the induced upper and lower subtractions are translation invariant with respect to~$\oplus$,
which permits us to introduce upper tolerances for elements in an optimal solution
to \eqref{e:DO} and lower tolerances for elements outside the optimal solution and
establish precise formulas for their evaluation. They are exact counterparts for all
strict {\sf A}-operations of the corresponding formulae from \cite{Libura} established
in the simplest case of~$+$. The tolerance function for problem \eqref{e:DO} is
defined on $X$ via any optimal solution $S^*$ to \eqref{e:DO} as follows: its value
at $x\in S^*$ is equal to the upper tolerance of $x$ and its value at
$y\in X\setminus S^*$ is equal to the lower tolerance of~$y$. We prove that the
tolerance function is independent of a particular optimal solution to \eqref{e:DO},
and so, it is an invariant of the problem \eqref{e:DO} itself. Also, we show that the
tolerance function is useful in characterizing the whole set of optimal solutions and,
in particular, their uniqueness: the optimal solution to \eqref{e:DO} is unique iff the
tolerance function never vanishes.

It is to be noted that the case of nonstrict {\sf A}-operations $\oplus$ is more
intricate and is yet to be studied in detail. Following the tradition, we elaborate
only on the case of nonstrict {\sf A}-operation of addition $\oplus=\max$.
Also, we do not touch upon the computational complexity of the
sensitivity analysis in our approach, which is a separate problem in itself
(concerning the latter see \cite{Chakra,GS,Lawler,Libura10,Pardalos,Rama}).

The paper is organized as follows. In Section~\ref{s:aop} we develop a theory of
{\sf A}-ope\-rations and their two subtractions, upper and lower, to be applied
throughout the paper. In Section~\ref{s:DOP} we formulate a general optimization
problem for objective functions generated by {\sf A}-operations. In Section~\ref{s:usi}
we determine upper stability intervals and in Section~\ref{s:lt}---lower stability intervals of
costs, for which an optimal solution preserves its optimality. In Section~\ref{s:TF} for strict
{\sf A}-operations we define an invariant of the optimization problem---the tolerance
function, and establish its basic properties.

\section{{\sf A}-operations on $\Rb^+$ and their inverses.} \label{s:aop}

Throughout the paper $\Rb^+$ denotes the set $[0,\infty)$ of all nonnegative real
numbers or the set $(0,\infty)$ of all positive real numbers.

The aim of this section is to introduce an operation of generalized addition or generalized
multiplication on $\Rb^+$, called an {\sf A}-operation, and to gather its properties
for future reference.

\subsection{Definition.} \label{d:A}
A continuous function $A:[0,\infty)\times[0,\infty)\to[0,\infty)$ is said to be
an {\em{\sf A}-op\-er\-a\-tion\/} on $\Rb^+$ if $A$ maps $\Rb^+\times\Rb^+$
into $\Rb^+$ and, given $u,v,w\in[0,\infty)$, the following four conditions (axioms)
are satisfied:
  \begin{itemize}
\item[(A.1)] $A(u,A(v,w))=A(A(u,v),w)$ (associativity of $A$);
\item[(A.2)] $A(u,v)=A(v,u)$ (commutativity of $A$);
\item[(A.3)] $u<v$ implies $A(u,w)\le A(v,w)$ (monotonicity of $A$);
\item[(A.4)] $A(u,v)\to\infty$ as $u\to\infty$ for all $v\in\Rb^+$ (unboundedness of $A$).
  \end{itemize}
If, instead of (A.3), the function $A$ satisfies
  \begin{itemize}
\item[(A.3$_{\mbox{\footnotesize\rm s}}$)]
given $u,v\in[0,\infty)$ and $w\in\Rb^+$, $u<v$ implies $A(u,w)<A(v,w)$,
  \end{itemize}
then the {\sf A}-operation $A$ is said to be {\em strict\/} on~$\Rb^+$.

We denote by $\mathcal{A}(\Rb^+)$ the set of all {\sf A}-operations on $\Rb^+$ and
by $\mathcal{A}_{\mbox{\footnotesize\rm s}}(\Rb^+)$---the set of those operations
$A\in\mathcal{A}(\Rb^+)$, which are strict on~$\Rb^+$.

\subsection{{\!\!}} In extending an $A\in\mathcal{A}(\Rb^+)$ to any finite number of terms
it is convenient to set $u\oplus v\equiv u\oplus_Av=A(u,v)$ for all $u,v\in[0,\infty)$ and,
given $u_1,\dots,u_n\in[0,\infty)$ with an $n\in\Nb$, we put
  \begin{equation} \label{e:e12}
u\oplus\varnothing=u,\qquad\oplus_{i=1}^1 u_i=u_1,\qquad
\oplus_{i=1}^2u_i=u_1\oplus u_2,
  \end{equation}
and, inductively,
  \begin{equation} \label{e:e3}
\oplus_{i=1}^nu_i=A\bigl(\oplus_{i=1}^{n-1}u_i,u_n\bigr)=
\bigl(\oplus_{i=1}^{n-1}u_i\bigr)\oplus u_n\quad\mbox{if}\quad n\ge3.
  \end{equation}
We will write $A(u,v)$ or $u\oplus v$ indifferently as well as
$A=\oplus\in\mathcal{A}(\Rb^+)$. Now, conditions (A.1)--(A.4) can be rewritten
more commonly as
  \begin{eqnarray*}
&u\oplus(v\oplus w)=(u\oplus v)\oplus w,\quad u\oplus v=v\oplus u,\quad
u<v\implies u\oplus w\le v\oplus w,&\\[2pt]
&\mbox{and $u\oplus v\to\infty$ as $u\to\infty$
for all $v\in\Rb^+$, respectively.}&
  \end{eqnarray*}

\subsection{Generalized addition and multiplication.} \label{ss:gam}
Of particular importance on $\Rb^+=[0,\infty)$ are {\sf A}-operations $\oplus$,
generalizing the usual addition, which, along with axioms (A.1)--(A.3), satisfy the condition
  \begin{itemize}
\item[(A.5)] $0\oplus v=v$ for all $v\in[0,\infty)$.
  \end{itemize}
Note that the zero $0$ plays the role of the {\em neutral\/} element with respect to~$\oplus$,
axiom~(A.4) is redundant in this case (in fact, (A.3), (A.2) and (A.5) imply
$u\oplus v\ge u\oplus0=u$) and $u\oplus v\ge\max\{u,v\}$ for all $u,v\in[0,\infty)$.
{\sf A}-operations on $\Rb^+$, satisfying (A.5), are called $F$-operations in
\cite{MP74} (see also \cite{DAN06}, \cite{Folia08}, \cite[Section~3]{Musielak}).

The case $\Rb^+=(0,\infty)$ is more appropriate for {\sf A}-operations $\oplus$,
generalizing the usual multiplication, i.e., satisfying (A.1)--(A.4) and
   \begin{itemize}
\item[(A.6)] $0\oplus v=0$ for all $v\in(0,\infty)$.
  \end{itemize}
Also, we assume the existence of the {\em neutral\/} element (the unit with respect
to~$\oplus$) ${\pmb e}\in\Rb^+$ such that ${\pmb e}\oplus v=v$ for all $v\in[0,\infty)$.

\subsection{Convention.} \label{ss:conv}
In what follows {\sf A}-operations $A=\oplus$ on $\Rb^+=[0,\infty)$ or $(0,\infty)$
will be treated in a unified way, however, on $[0,\infty)$ they are assumed to satisfy
axioms (A.1)--(A.3) and (A.5) and are called {\em{\sf A}-operations of addition\/}
and on $(0,\infty)$---axioms (A.1)--(A.4) and (A.6) and are called
{\em{\sf A}-operations of multiplication}.

\subsection{Examples of {\sf A}-operations.} \label{ss:exA}
Let $u,v\in[0,\infty)$ and $p\in(0,\infty)$. The following are examples of
{\sf A}-operations of addition on $\Rb^+=[0,\infty)$
(cf. \cite{Folia08}, \cite[Chapter~3]{Musielak}):
  \begin{itemize}
\item[(a$_1$)] $A_1(u,v)=u+v$ (the usual addition operation);
\item[(a$_2$)] $A_2(u,v)=(u^p+v^p)^{1/p}$;
\item[(a$_3$)] $A_3(u,v)=\max\{u,v\}=\frac12(u+v+|u-v|)$;
\item[(a$_4$)] $A_4(u,v)=\frac1p\log(e^{pu}+e^{pv}-1)$;
\item[(a$_5$)] $A_5(u,v)=u+v+puv$;
\item[(a$_6$)] $A_6(u,v)=G(\max\{u,v\},u+v)$, where $G(a,b)$ is defined
  for $0\le a\le b$ by: $G(a,b)=b$ if $b<1$, $G(a,b)=1$ if $a<1$ and $b\ge1$,
  and $G(a,b)=a$ if $a\ge1$.
  \end{itemize}

Examples of {\sf A}-operations of multiplication on $\Rb^+=(0,\infty)$ are as follows:
  \begin{itemize}
\item[(a$_7$)] $A_7(u,v)=uv$ (the usual operation of multiplication)
  with the unit ${\pmb e}=1$;
\item[(a$_8$)] $A_8(u,v)=puv$, the neutral element being ${\pmb e}=1/p$;
\item[(a$_9$)] $A_9(u,v)=\frac1p\log\bigl(1+(e^{pu}-1)(e^{pv}-1)\bigr)$
  with the unit ${\pmb e}=(\log2)/p$;
\item[(a$_{10}$)] $A_{10}(u,v)=\frac1p\bigl(\exp[\log(1+pu)\cdot\log(1+pv)]-1\bigr)$,
  the unit being ${\pmb e}=(e-1)/p$.
  \end{itemize}

That all these ten operations are indeed {\sf A}-operations on $\Rb^+$ (satisfying
convention~\ref{ss:conv}) will be more clear from Section~\ref{ss:equiv}.

Operations $A_1$, $A_2$, $A_4$ and $A_5$ are strict on $[0,\infty)$, operations
$A_7$, $A_8$, $A_9$ and $A_{10}$ are strict on $(0,\infty)$, while operations
$A_3=\max$ and $A_6$ are {\em not\/} strict on $[0,\infty)$.

\smallbreak
For several elements $u_1,\dots,u_n\in[0,\infty)$ operations $A_1$ through $A_{10}$
extend in the way exposed in \eqref{e:e3} and Sections~\ref{ss:equiv} and
\ref{ss:eqvD} (cf.\ expression for $\widehat\oplus_{i=1}^nu_i$); see also
Examples~\ref{ss:eom}.

\subsection{Two properties.} \label{ss:twop}
Two simple properties of {\sf A}-operations on $\Rb^+$ are straightforward: if
$\oplus\in\mathcal{A}(\Rb^+)$ and $u,u_1,v,v_1\in[0,\infty)$, then
  \begin{equation*}
\mbox{$u\le u_1$ \,and \,$v\le v_1$ \,imply \,$u\oplus v\le u_1\oplus v_1$,}
  \end{equation*}
and if, in addition, $\oplus$ is strict on $\Rb^+$, then, given $w\in\Rb^+$,
  \begin{equation} \label{e:canl}
\mbox{$u\oplus w\le v\oplus w$ \,implies \,$u\le v$ \,(cancellation law),}
  \end{equation}
which, in particular, gives: if $u\oplus w=v\oplus w$, then $u=v$.

\subsection{Generating {\sf A}-operations.} \label{ss:equiv}
Following \cite{Folia08} or \cite[Section~3]{Musielak}, here we introduce an equivalence
relation on the set $\mathcal{A}(\Rb^+)$.

Let $\varphi:[0,\infty)\to[0,\infty)$ be a continuous strictly increasing function vanishing
at zero (only) and such that $\varphi(u)\to\infty$ as $u\to\infty$ (such functions are
said to be {\em$\varphi$-functions}, cf.\ \cite{Folia08}, \cite{Musielak}).
Given $A\in\mathcal{A}(\Rb^+)$ and $u,v\in\Rb^+$, we set
  \begin{equation*}
(E_\varphi(A))(u,v)=\varphi^{-1}\bigl(A(\varphi(u),\varphi(v))\bigr),
  \end{equation*}
where $\varphi^{-1}:[0,\infty)\to[0,\infty)$ is the inverse function of $\varphi$.
Clearly, $E_\varphi(A)\in\mathcal{A}(\Rb^+)$, and so, $E_\varphi$ maps
$\mathcal{A}(\Rb^+)$ into itself.

If $\mbox{\rm id}_X$ denotes the identity map of a set $X$ (i.e.,
$\mbox{\rm id}_X(x)=x$ for all $x\in X$) and $\varphi_0=\mbox{\rm id}_{[0,\infty)}$,
then $E_{\varphi_0}=\mbox{\rm id}_{\mathcal{A}(\Rb^+)}$. Also, given two
$\varphi$-functions $\varphi$ and $\psi$, we find
$E_{\varphi\circ\psi}=E_\varphi\circ E_\psi$, where $\circ$ designates the usual
composition of maps. It follows that the relation $\sim$ on $\mathcal{A}(\Rb^+)$
defined for $A,B\in\mathcal{A}(\Rb^+)$ by
  \begin{equation*}
\mbox{$B\sim A$ \,\,iff \,\,there exists a $\varphi$-function $\varphi$ such that $B=E_\varphi(A)$},
  \end{equation*}
is an equivalence relation on $\mathcal{A}(\Rb^+)$. The equivalence class $[A]$ of an
{\sf A}-operation $A$ under $\sim$ is given by
$[A]=\{E_\varphi(A):\mbox{$\varphi$ is a $\varphi$-function}\}$. It is to be noted that
if $A\in\mathcal{A}_{\mbox{\footnotesize\rm s}}(\Rb^+)$, then
$[A]\subset\mathcal{A}_{\mbox{\footnotesize\rm s}}(\Rb^+)$.

Now we turn back to Examples~\ref{ss:exA}. We have:
  \begin{eqnarray*}
&&\mbox{$A_2=E_\varphi(A_1)$ \,and \,$A_7=E_\varphi(A_7)$ \,\,with
  \,$\varphi(u)=u^p$,}\\[2pt]
&&\mbox{$A_4=E_\psi(A_1)$ \,and \,$A_9=E_\psi(A_7)$ \,\,with
  \,$\psi(u)=e^{pu}-1$,}\\[2pt]
&&\mbox{$A_5=E_\chi(A_1)$ \,and \,$A_{10}=E_\chi(A_7)$ \,with
  \,$\chi(u)=\log(1+pu)$,}
  \end{eqnarray*}
and $A_8=E_\xi(A_7)$ with $\xi(u)=pu$, where $p>0$ and $u\in[0,\infty)$.
Thus, $A_1\sim A_2\sim A_4\sim A_5$ and $A_7\sim A_8\sim A_9\sim A_{10}$,
while $A_1$, $A_3$, $A_6$ and $A_7$ are not mutually equivalent under~$\sim$.
Note also that the equivalence class of $A_3=\max$ is $[A_3]=\{A_3\}$.

If ${\pmb e}\in\Rb^+$ is the neutral element with respect to the {\sf A}-operation
$\oplus$ and $\widehat\oplus=E_\varphi(\oplus)$ for a $\varphi$-function $\varphi$,
then $\widehat{\pmb e}=\varphi^{-1}({\pmb e})$ is the neutral element with respect
to the {\sf A}-operation~$\widehat{\oplus}$: in fact, given $u\ge0$, we find
  \begin{equation*}
\widehat{\pmb e}\widehat{\oplus}u
=\varphi^{-1}\bigl(\varphi(\widehat{\pmb e})\oplus\varphi(u)\bigr)
=\varphi^{-1}\bigl({\pmb e}\oplus\varphi(u)\bigr)=\varphi^{-1}(\varphi(u))=u.
  \end{equation*}

\subsection{Upper and lower subtractions.} \label{ss:uplo}
Here we treat two inverses of an {\sf A}-operation $\oplus$, also called generalized
subtractions (or divisions). Having an equation of the form $u\oplus v=w$ (or
inequality $u\oplus v\le w$) with $u,v,w\in\Rb^+$, we are going to write
$u=w\ominus v$ (or $u\le w\ominus v$), so that we would get
$(w\ominus v)\oplus v=w$ (or $(w\ominus v)\oplus v\le w$, respectively).
We will achieve this by introducing two ``operations'' of upper and lower subtractions
on $\Rb^+$, $\upsub$ and $\losub$, as follows.

Suppose the {\sf A}-operation $A=\oplus\in\mathcal{A}(\Rb^+)$ is given.

We define the domain $D(\upsub)$ of the {\em upper subtraction\/} $\upsub$
for $\oplus$ by
  \begin{equation*}
D(\upsub)=\bigl\{(w,v)\in\Rb^+\times\Rb^+:\mbox{$\exists\,u_0\in\Rb^+$
such that $u_0\oplus v\le w$}\bigr\},
  \end{equation*}
and, given $(w,v)\in D(\upsub)$, the value $w\upsub v\in\Rb^+$ is defined by
  \begin{equation} \label{e:dups}
w\upsub v=\max\{u\in\Rb^+:u\oplus v\le w\}=\max\{u\in\Rb^+:u\oplus v=w\}.
  \end{equation}

The {\em lower subtraction\/} $\losub$ for $\oplus$ is defined for all
$(w,v)\in\Rb^+\times\Rb^+$ by
  \begin{equation} \label{e:dlos}
w\losub v=\min\{u\in\Rb^+:u\oplus v\ge w\}\in\Rb^+.
  \end{equation}

Definitions \eqref{e:dups} and \eqref{e:dlos} are similar to taking the right or left inverse
of a not necessarily strictly increasing function on $[0,\infty)$ (cf.\ \cite{KR}),
depending on a parameter.

Our primary aim now is to verify that these definitions are correct.

\subsection{Lemma.} \label{l:ucorr} 
{\em The subtractions\/ $\upsub$ and\/ $\losub$ are well defined.}
\medbreak
\proof{Proof.}
1. Let us show that definition \eqref{e:dups} is correct. First, we note that the domain
$D(\upsub)$ is non\-empty. In fact, given $u,v\in\Rb^+$, setting $u_0=u$ and
$w=u\oplus v$, we find $u_0\oplus v=w$, and so, $(u\oplus v,v)=(w,v)\in D(\upsub)$.

Now, let $(w,v)\in D(\upsub)$, and so, there exists $u_0\in\Rb^+$ such that
$u_0\oplus v\le w$. Defining the set $\overline U=\{u\in\Rb^+:u\oplus v\le w\}$, we find
$\overline U\ne\varnothing$ (since $u_0\in\overline U$), $\overline U$ is bounded
from below (because $\overline U\subset\Rb^+$) and bounded from above
(by axiom (A.4)) and $\overline U$ is closed in $\Rb$ (by the continuity of the
function $[u\mapsto u\oplus v]:\Rb^+\to\Rb^+$), and so, $\overline U$ is compact
in~$\Rb$. We set $\overline u=\sup\overline U=\max\overline U\in\Rb$ and assert
that $\overline u\in\Rb^+$ and $\overline u\oplus v=w$. In fact, given
$u\in\overline U$, we have $u\le\overline u$ and, since $u_0\in\overline U$, inequality
$u_0\le\overline u$ implies $\overline u\in\Rb^+$. Since $\overline u=\max\overline U$,
we have $\overline u\oplus v\le w$. If we assume that $\overline u\oplus v<w$, then
we note that, by (A.4), there exists $\widehat u\in\Rb^+$ such that $w<\widehat u\oplus v$
and, moreover, by (A.3), $\overline u<\widehat u$; now, by the intermediate value theorem,
there exists $\widetilde u\in\Rb^+$ with $\overline u<\widetilde u<\widehat u$ such that
$\widetilde u\oplus v=w$, and so, $\widetilde u\in\overline U$ and
$\max\overline U=\overline u<\widetilde u$, which is a contradiction. It follows that
$\overline u\oplus v=w$ and, hence, $w\upsub v=\overline u$.

2. Let us show that the quantity \eqref{e:dlos} is well defined. Given
$(w,v)\in\Rb^+\times\Rb^+$, we define the set
$\underline U=\{u\in\Rb^+:u\oplus v\ge w\}$.
Then $\underline U\ne\varnothing$ (by virtue of (A.4)) and $\underline U$ is bounded
from below (since $\underline U\subset\Rb^+$), and so,
$\underline u=\inf\underline U\in[0,\infty)=\{0\}\cup\Rb^+$.
If $\Rb^+=[0,\infty)$, then we note that, by the continuity of the function
$u\mapsto u\oplus v$, the set $\underline U$ is closed in $[0,\infty)$, and so,
$\underline u=\min\underline U$ and $\underline u\in[0,\infty)=\Rb^+$.
In the case when $\Rb^+=(0,\infty)$, we have
  \begin{equation} \label{e:ugeu}
u\oplus v\ge w\quad\mbox{for \,all}\quad u>\underline u,
  \end{equation}
and, by virtue of (A.6), $u\oplus v\to0\oplus v=0$ as $u\to0$, and so, there
exists $u_1>0$, depending on $w>0$, such that $u_1\oplus v<w$. It follows from
\eqref{e:ugeu} that $\underline u\ge u_1$, i.e., $\underline u>0$. Passing to the limit
as $u\to\underline u$ in \eqref{e:ugeu}, by the continuity of $u\mapsto u\oplus v$,
we get $\underline u\oplus v\ge w$, and so, $\underline u\in\underline U$ and
$\underline u=\min\underline U\in(0,\infty)=\Rb^+$. Thus, $w\losub v=\underline u$.
\Halmos\endproof

\subsection{Remark.} \label{r:rev}
We have shown in step 2 in the above proof that if $\Rb^+=(0,\infty)$, then (A.6) implies
$\inf\underline U=\min\underline U$ for all $w>0$. It is to be noted that the converse
implication holds as well. In fact, given $\vep>0$, we set
$u_\vep=\inf\underline U=\min\underline U$, where
$\underline U=\{u\in(0,\infty):u\oplus v\ge\vep\}$. We assert that $u\oplus v<\vep$
for all $0<u<u_\vep$ (and so, $0\oplus v=\lim_{u\to0}u\oplus v=0$ implying (A.6)):
indeed, if we assume that $u\oplus v\ge\vep$, then $u\in\underline U$, and so,
$u\ge\min\underline U=u_\vep$, which contradicts the inequality $u<u_\vep$.

\smallbreak
Note that the value $w\upsub v$, as opposed to $w\losub v$, may {\em not\/} be
defined for all $w,v\in\Rb^+$. A comparison of the two subtractions $\upsub$ and
$\losub$ for a $\oplus\in\mathcal{A}(\Rb^+)$ is given in the following

\subsection{Lemma.} \label{l:comp} {\em
{\rm(a)} If $(w,v)\in D(\upsub)$, then $w\losub v\le w\upsub v$.
  \begin{itemize}
\item[{\rm(b)}] If\/ $\oplus$ is an\/ {\sf A}-operation of addition on\/ 
  $\Rb^+=[0,\infty)$, then $(w,v)\in D(\upsub)$ iff\/ $0\le v\le w$,
  and inequalities $v\ge w\ge0$ imply $w\losub v=0$.
\item[{\rm(c)}] If\/ $\oplus$ is an\/ {\sf A}-operation of multiplication on\/ 
  $\Rb^+=(0,\infty)$, then $(w,v)\in D(\upsub)$ for all $w,v>0$.
  \end{itemize}
}
\medbreak
\proof{Proof.}
(a) is a straightforward consequence of \eqref{e:dups} and \eqref{e:dlos}.
\smallbreak
(b) Given $(w,v)\in D(\upsub)$, there exists $u_0\ge0$ such that $u_0\oplus v\le w$,
and so, by virtue of (A.5) and (A.3), $v=0\oplus v\le u_0\oplus v\le w$. On the other
hand, if $v\le w$, then, by (A.5), $0\oplus v=v\le w$, and so, $(w,v)\in D(\upsub)$.

Suppose $v\ge w$. By (A.3) and (A.5), $u\oplus v\ge u\oplus w\ge0\oplus w=w$ for
all $u\ge0$, and so, $\underline U=\{u\ge0:u\oplus v\ge w\}=[0,\infty)$ and
$w\losub v=\min\underline U=0$.
\smallbreak
(c) Given $w,v>0$, by virtue of (A.6), $u\oplus v\to0$ as $u\to0$, and so,
there exists $u_0>0$ such that $u_0\oplus v<w$ implying $(w,v)\in D(\upsub)$.
\Halmos\endproof

\subsection{Examples of upper and lower subtractions.} \label{ss:exul}
Most examples in this paper will be demonstrated for the three basic (representatives
of equivalence classes of) {\sf A}-operations $A_2$, $A_3$ and $A_7$ from
Section~\ref{ss:exA}, i.e., $u\oplus v=(u^p+v^p)^{1/p}$ with $p>0$,
$u\oplus v=\max\{u,v\}$ and $u\oplus v=uv$, where $u,v\ge0$, also abbreviated
as {\em $p$-sum, max\/} and {\em product\/} operations, respectively.

Let $\oplus\in\mathcal{A}(\Rb^+)$ and $\upsub$ and $\losub$ be the upper and lower
subtractions for~$\oplus$.

(a) Suppose $u\oplus v=(u^p+v^p)^{1/p}$ on $\Rb^+=[0,\infty)$ with $p>0$. We have:
  \begin{equation*}
\mbox{$(w,v)\in D(\upsub)$ \,iff \,$0\le v\le w$, \,and \,$w\upsub v=(w^p-v^p)^{1/p}$;}
  \end{equation*}
  \begin{equation*}
\mbox{given \,$w,v\ge0$,}\quad w\losub v=\left\{
  \begin{array}{ccl}
(w^p-v^p)^{1/p} & \mbox{if} & v<w,\\
0                       & \mbox{if} & v\ge w.
  \end{array}\right.
  \end{equation*}
Note that $w\upsub v=w\losub v$ for all $0\le v\le w$.
\smallbreak
(b) Let $u\oplus v=\max\{u,v\}$ on $\Rb^+=[0,\infty)$. Then
  \begin{equation*}
\mbox{$(w,v)\in D(\upsub)$ \,iff \,$0\le v\le w$, \,and \,$w\upsub v=w$;}
  \end{equation*}
  \begin{equation*}
\mbox{given \,$w,v\ge0$,}\quad w\losub v=\left\{
  \begin{array}{ccl}
w & \mbox{if} & v<w,\\
0 & \mbox{if} & v\ge w.
  \end{array}\right.
  \end{equation*}
In fact, if $v\le w$, then $\overline U=\{u\ge0:\max\{u,v\}\le w\}=[0,w]$, and so,
by \eqref{e:dups}, $w\upsub v=\max\overline U=w$. By virtue of
Lemma~\ref{l:comp}(b), $v\ge w$ implies $w\losub v=0$, and if $v<w$, then
$u\in\underline U$ iff $\max\{u,v\}\ge w$ iff $u\ge w$, and so,
$\underline U=[w,\infty)$ and $w\losub v=\min\underline U=w$.

Note that $w\upsub v=w=w\losub v$ for all $0\le v<w$ and $0\upsub0=0=0\losub0$,
while if $w>0$, then we have $w\losub w=0<w=w\upsub w$.

(c) Suppose $u\oplus v=uv$ on $\Rb^+=(0,\infty)$. It follows from
Lemma~\ref{l:comp}(c) that $D(\upsub)=(0,\infty)\times(0,\infty)$ and
$w\upsub v=w\losub v=w/v$ is the usual operation of division for all $w,v>0$.

\smallbreak
Now we establish basic (in)equalities related to the upper and lower subtractions
$\upsub$ and $\losub$ for an {\sf A}-operation~$\oplus$.

\subsection{Lemma.} \label{l:ineq} {\em
Given\/ $\oplus\in\mathcal{A}(\Rb^+)$, we have\/{\rm:}
  \begin{eqnarray}
(w\upsub v)\oplus v&=&w \qquad\qquad\qquad\quad
  \quad\mbox{for \,all \,\,\,$(w,v)\in D(\upsub);$} \label{e:ups}\\[1pt]
&&w\,\le\,(w\losub v)\oplus v
  \qquad\mbox{for \,all \,\,\,$w,v\in\Rb^+;$} \label{e:los}\\[2pt]
(w\oplus v)\losub v&\le&w\,\le\,(w\oplus v)\upsub v
  \qquad\mbox{for \,all \,\,\,$w,v\in\Rb^+.$} \label{e:wvv}
  \end{eqnarray}
Also, the following criterion holds\/{\rm:}
  \begin{equation} \label{e:u=l}
\mbox{$\oplus$ is strict on $\Rb^+$ \,iff \,\,$w\upsub v=w\losub v$ \,for \,all
\,$(w,v)\in D(\upsub)$.}
  \end{equation}
}
\proof{Proof.}
1. Equality \eqref{e:ups} is the characterizing property of the quantity $w\upsub v$,
which follows immediately from \eqref{e:dups}.

2. Inequality \eqref{e:los} is the characterizing property of $w\losub v$, which is a
consequence of \eqref{e:dlos} (cf.\ also Remark~\ref{r:214}(a)).

3. The left-hand side inequality in \eqref{e:wvv} follows from \eqref{e:dlos}:
  \begin{equation*}
(w\oplus v)\losub v=\min\{u\in\Rb^+:u\oplus v\ge w\oplus v\}\le w.
  \end{equation*}
Setting $u_0=w$, we find $u_0\oplus v=w\oplus v$, and so, $(w\oplus v,v)\in D(\upsub)$,
and the right-hand side inequality in \eqref{e:wvv} follows from \eqref{e:dups}:
  \begin{equation*}
(w\oplus v)\upsub v=\max\{u\in\Rb^+:u\oplus v=w\oplus v\}\ge w.
  \end{equation*}

4. Let us prove \eqref{e:u=l}. Let $\oplus$ be strict. Given $(w,v)\in D(\upsub)$, by
virtue of Lemma~\ref{l:comp}(a), we have to show only that $w\losub v\ge w\upsub v$.
In fact, we assert that $u\ge w\upsub v$ for all
$u\in\underline U=\{u\in\Rb^+:u\oplus v\ge w\}$, for, otherwise, if $u<w\upsub v$,
then (A.3$_{\mbox{\footnotesize\rm s}}$) and \eqref{e:ups} imply
$u\oplus v<(w\upsub v)\oplus v=w$, which is a contradiction. Thus,
$w\losub v=\min\underline U\ge w\upsub v$.

Suppose the equality on the right in \eqref{e:u=l} holds, and let $u>v\ge0$ and
$w\in\Rb^+$. It follows from (A.3) that $u\oplus w\ge v\oplus w$. If we assume that
$u\oplus w=v\oplus w$, then, by \eqref{e:dups} and the left-hand side inequality
in \eqref{e:wvv}, we get
  \begin{equation*}
u\le(v\oplus w)\upsub w=(v\oplus w)\losub w\le v,
  \end{equation*}
which contradicts to $u>v$. Thus, $u\oplus w>v\oplus w$, and
(A.3$_{\mbox{\footnotesize\rm s}}$) follows.
\Halmos\endproof

\subsection{Remarks.} \label{r:214}
(a) Strict inequality may hold in \eqref{e:los} even for strict {\sf A}-operations:
if $\oplus$ is as in Lemma~\ref{l:comp}(b) and $v>w\ge0$, then
  \begin{equation*}
w<v=0\oplus v=(w\losub v)\oplus v.
  \end{equation*}

Also, one cannot replace the inequality $\ge w$ in \eqref{e:dlos} by the equality $=w$:
in fact, if $\oplus=\max$, then taking into account Example~\ref{ss:exul}(b),
we have, for $v>w\ge0$,
  \begin{equation*}
w\losub v=\min\{u\ge0:\max\{u,v\}\ge w\}=0,
  \end{equation*}
whereas $\{u\ge0:\max\{u,v\}=w\}=\varnothing$.

(b) If $\oplus$ is strict on $\Rb^+$, then, by \eqref{e:u=l}, equalities hold in \eqref{e:wvv}.
However, if $\oplus$ is not strict, then inequalities in \eqref{e:wvv} may be strict. By
virtue of Example~\ref{ss:exul}(b), for $\oplus=\max$ we have: if $v=w>0$, then
  \begin{equation*}
(w\oplus v)\losub v=(w\oplus w)\losub w=\max\{w,w\}\losub w=w\losub w=0<w,
  \end{equation*}
and if $v>w\ge0$, then
  \begin{equation*}
(w\oplus v)\upsub v=\max\{w,v\}\upsub v=v\upsub v=v>w.
  \end{equation*}
In this respect we note that, by Example~\ref{ss:exul}(b), $w\upsub v=w$ for all
$0\le v\le w$, and so, equality \eqref{e:ups} is of the form
  \begin{equation*}
(w\upsub v)\oplus v=w\oplus v=\max\{w,v\}=w.
  \end{equation*}

The following lemma shows that the functions $w\mapsto w\upsub v$ and
$w\mapsto w\losub v$ are nondecreasing (in the first variable).

\subsection{Lemma.} \label{l:fivar} {\em
Let\/ $\oplus\in\mathcal{A}(\Rb^+)$, $w_1,w_2,v\in\Rb^+$ and $w_1\le w_2$.
We have\/{\rm:}
  \begin{itemize}
\item[{\rm(a)}] if $(w_1,v)\in D(\upsub)$, then $(w_2,v)\in D(\upsub)$ and
  $w_1\upsub v\le w_2\upsub v;$
\item[{\rm(b)}] $w_1\losub v\le w_2\losub v$.
  \end{itemize}
}
\medbreak
\proof{Proof.}
(a) Condition $(w_1,v)\in D(\upsub)$ implies the existence of $u_0\in\Rb^+$ such that
$u_0\oplus v\le w_1$, which gives $u_0\oplus v\le w_2$, and so, $(w_2,v)\in D(\upsub)$.
Setting $u_1=w_1\upsub v$, by virtue of \eqref{e:ups}, we find
  \begin{equation*}
u_1\oplus v=(w_1\upsub v)\oplus v=w_1\le w_2,
  \end{equation*}
and so, \eqref{e:dups} yields $u_1\le w_2\upsub v$.
\smallbreak
(b) If $u_2=w_2\losub v$, then it follows from \eqref{e:los} that
  \begin{equation*}
w_1\le w_2\le(w_2\losub v)\oplus v=u_2\oplus v,
  \end{equation*}
whence \eqref{e:dlos} implies $w_1\losub v\le u_2$.
\Halmos\endproof

Several more inequalities will be needed in the sequel. If $u$, $v$ and $w$ are real numbers,
then for the usual operations of addition $+$ and subtraction $-$ we have: $u\le v$
implies $w-v\le w-u$, and $w-(v-u)=(w+u)-v$. In the next two lemmas we study to what
extent these two properties carry over to general {\sf A}-operations
$\oplus\in\mathcal{A}(\Rb^+)$ and their upper and lower subtractions $\upsub$
and~$\losub$.

\subsection{Lemma.} \label{l:scnde} {\em
Suppose $u,v,w\in\Rb^+$ and $u\le v$. Then we have\/{\rm:}
  \begin{itemize}
\item[{\rm(a)}] if $(w,v)\in D(\upsub)$, then $(w,u)\in D(\upsub)$ and
  $w\upsub v\le w\upsub u;$
\item[{\rm(b)}] $w\losub v\le w\losub u$.
  \end{itemize}
}
\medbreak
\proof{Proof.}
(a) By virtue of (A.3) and \eqref{e:ups}, we get
  \begin{equation*}
(w\upsub v)\oplus u\le(w\upsub v)\oplus v=w,
  \end{equation*}
and so, $(w,u)\in D(\upsub)$, and the definition \eqref{e:dups} of $w\upsub u$
(i.e., the maximality of $w\upsub u$) implies the desired inequality in~(a).

(b) It follows form (A.3) and \eqref{e:los} that
  \begin{equation*}
(w\losub u)\oplus v\ge(w\losub u)\oplus u\ge w,
  \end{equation*}
and the minimality of $w\losub v$ from \eqref{e:dlos} gives the inequality in~(b).
\Halmos\endproof

\subsection{Lemma.} \label{l:luineq} {\em
Let $\oplus\in\mathcal{A}(\Rb^+)$ and $u,v,w\in\Rb^+$. We have\/{\rm:}
  \begin{itemize}
\item[{\rm(a)}] if\/ $(v,u)\in D(\upsub)$ and $(w,v\upsub u)\in D(\upsub)$, then
$(w\oplus u,v)\in D(\upsub)$ and
  \begin{equation*}
(w\oplus u)\losub v\le w\upsub(v\upsub u)\le(w\oplus u)\upsub v;
  \end{equation*}
\item[{\rm(b)}] if\/ $(v,u)\in D(\upsub)$, then
  $(w\oplus u)\losub v\le w\losub(v\upsub u);$
\item[{\rm(c)}] if\/ $(w,v\losub u)\in D(\upsub)$, then $(w\oplus u,v)\in D(\upsub)$
  and $w\upsub(v\losub u)\le(w\oplus u)\upsub v;$
\item[{\rm(d)}] $(w\oplus u)\losub[(v\losub u)\oplus u]\le w\losub(v\losub u)$.
  \end{itemize}
}
\medbreak
\proof{Proof.}
(a) Setting $w_1=w\upsub(v\upsub u)$, by virtue of \eqref{e:ups}, we find
  \begin{equation*}
w=[w\upsub(v\upsub u)]\oplus(v\upsub u)=w_1\oplus(v\upsub u),
  \end{equation*}
and so, once again \eqref{e:ups} implies
  \begin{equation*}
w\oplus u=w_1\oplus(v\upsub u)\oplus u=w_1\oplus v.
  \end{equation*}
Consequently, $(w\oplus u,v)\in D(\upsub)$. Applying \eqref{e:wvv}, we get
  \begin{equation*}
(w\oplus u)\losub v=(w_1\oplus v)\losub v\le w_1\le(w_1\oplus v)\upsub v
=(w\oplus u)\upsub v.
  \end{equation*}

(b) Set $w_2=w\losub(v\upsub u)$. It follows from \eqref{e:los} that
  \begin{equation*}
w\le[w\losub(v\upsub u)]\oplus(v\upsub u)=w_2\oplus(v\upsub u),
  \end{equation*}
and so, by (A.3) and \eqref{e:ups},
  \begin{equation*}
w\oplus u\le w_2\oplus(v\upsub u)\oplus u=w_2\oplus v.
  \end{equation*}
The desired inequality $(w\oplus u)\losub v\le w_2$ follows from definition~\eqref{e:dlos}.

\smallbreak
(c) It follows from \eqref{e:ups} that if $w_3=w\upsub(v\losub u)$, then
  \begin{equation*}
w=[w\upsub(v\losub u)]\oplus(v\losub u)=w_3\oplus(v\losub u),
  \end{equation*}
and so, \eqref{e:los} implies
  \begin{equation*}
w\oplus u=w_3\oplus(v\losub u)\oplus u\ge w_3\oplus v.
  \end{equation*}
Thus, $(w\oplus u,v)\in D(\upsub)$ and, by \eqref{e:dups}, $w_3\le(w\oplus u)\upsub v$.

\smallbreak
(d) Setting $w_4=w\losub(v\losub u)$ and applying \eqref{e:los}, we have
  \begin{equation*}
w\le[w\losub(v\losub u)]\oplus(v\losub u)=w_4\oplus(v\losub u),
  \end{equation*}
and so, by (A.3), $w\oplus u\le w_4\oplus(v\losub u)\oplus u$. Now, the desired inequality
in (d) follows from definition~\eqref{e:dlos}.
\Halmos\endproof

In the final lemma of this section we address the {\em translation invariance\/} of
subtractions $\upsub$ and $\losub$ with respect to the {\sf A}-operation~$\oplus$
that generates them.

\subsection{Lemma.} \label{l:trinv} {\em
Suppose $\oplus$ is an\/ {\sf A}-operation on\/ $\Rb^+$.
  \begin{itemize}
\item[{\rm(a)}] If $u\in\Rb^+$ and $(w,v)\in D(\upsub)$, then
$(u\oplus w,u\oplus v)\in D(\upsub)$ and
  \begin{equation} \label{e:oti}
w\upsub v\le(u\oplus w)\upsub(u\oplus v).
  \end{equation}
In addition, if $\oplus$ is strict, then we have the equality
$(u\oplus w)\upsub(u\oplus v)=w\upsub v$.
\item[{\rm(b)}] If $u,v,w\in\Rb^+$, then $(u\oplus w)\losub(u\oplus v)\le w\losub v$. If,
in addition, $\oplus$ is strict, then we have the equality
$(u\oplus w)\losub(u\oplus v)=w\losub v$.
  \end{itemize}
}
\medbreak
\proof{Proof.}
(a) Since $(w,v)\in D(\upsub)$, there exists $u_0\in\Rb^+$ such that $u_0\oplus v\le w$,
and so, by (A.1)--(A.3), $u_0\oplus(u\oplus v)=u\oplus(u_0\oplus v)\le u\oplus w$, which
implies that the pair $(u\oplus w,u\oplus v)$ is in the domain of~$\upsub$.

Taking into account \eqref{e:ups} and (A.2), we get
  \begin{equation*}
(w\upsub v)\oplus(u\oplus v)=[(w\upsub v)\oplus v]\oplus u=w\oplus u=u\oplus w,
  \end{equation*}
and so, inequality \eqref{e:oti} is a consequence of definition \eqref{e:dups}.

Now, suppose $\oplus$ is strict. Setting $u_1=(u\oplus w)\upsub(u\oplus v)$, by
\eqref{e:ups}, we have
  \begin{equation*}
u_1\oplus u\oplus v=[(u\oplus w)\upsub(u\oplus v)]\oplus(u\oplus v)=u\oplus w.
  \end{equation*}
By virtue of \eqref{e:canl}, we cancel by $u$ and get $u_1\oplus v=w$. Since, again by
\eqref{e:ups}, $(w\upsub v)\oplus v=w$, we find $u_1\oplus v=w=(w\upsub v)\oplus v$,
and so, cancelling by $v$, we arrive at $u_1=w\upsub v$, which is the desired equality.

\smallbreak
(b) Set $u_2=(u\oplus w)\losub(u\oplus v)$. It follows from (A.1)--(A.3) and \eqref{e:los}
that
  \begin{equation*}
(w\losub v)\oplus(u\oplus v)=[(w\losub v)\oplus v]\oplus u\ge w\oplus u=u\oplus w,
  \end{equation*}
and so, definition \eqref{e:dlos} implies $u_2=(u\oplus w)\losub(u\oplus v)\le w\losub v$.

Let $\oplus$ be strict. By virtue of \eqref{e:los}, we get
  \begin{equation*}
u_2\oplus u\oplus v=[(u\oplus w)\losub(u\oplus v)]\oplus(u\oplus v)\ge u\oplus w,
  \end{equation*}
and so, cancelling by $u$, we find $u_2\oplus v\ge w$, which, by virtue of definition
\eqref{e:dlos}, implies $w\losub v\le u_2$, and the desired equality readily follows.
\Halmos\endproof

\subsection{Remark.} \label{r:unstr}
The inequalities in Lemma~\ref{l:trinv} may be strict if $\oplus$ is not necessarily strict.
To see this, we set $\oplus=\max$ and take into account Example~\ref{ss:exul}(b):
given $0\le v<w<u$, we have
  \begin{equation*}
w\upsub v=w<u=u\upsub u=\max\{u,w\}\upsub\max\{u,v\}=(u\oplus w)\upsub(u\oplus v)
  \end{equation*}
and
  \begin{equation*}
(u\oplus w)\losub(u\oplus v)=\max\{u,w\}\losub\max\{u,v\}=u\losub u=0<w=w\losub v.
  \end{equation*}

\pagebreak
\section{Optimization problems.} \label{s:DOP}

In this section we introduce notation, definitions and assumptions to be used
throughout the paper.

\subsection{Optimization space.} \label{ss:os}
Let $X$ be a finite set of cardinality $|X|\ge2$ (e.g., $X=\{1,2,\dots,n\}$ with $n\ge2$),
called the {\em ground set}, $2^X$ be the family of all subsets of $X$ (i.e., the power
set of $X$) and $\dot2^X=2^X\setminus\{\varnothing\}$. For instance, a ground set
may be thought of as the collection of all edges of a graph (or arcs in the directed case).
Given a nonempty set $Y$, we denote, as usual, by $Y^X$ the set of all functions (maps)
$g:X\to Y$ from $X$ into~$Y$.

A {\em set of trajectories\/} on $X$ (or a canonical collection on $X$, cf.~\cite{JOGO12})
is a collection $\CS\subset\dot2^X$ of subsets of $X$ such that
  \begin{equation} \label{e:traj}
\cup\CS=X\,\quad\,\mbox{and}\,\quad\,\cap\CS=\varnothing,
  \end{equation}
where $\cup\CS$ is the union of $\CS$ ($=$\,the set of all $x\in X$ such that
$x\in S$ for some $S\in\CS$) and $\cap\CS$ is the intersection of $\CS$ ($=$\,the set
of all $x\in X$ such that $x\in S$ for all $S\in\CS$). It follows immediately from \eqref{e:traj}
that there are at least two trajectories in $\CS$, and so, $2\le|\CS|<2^{|X|}$.

A pair $(X,\CS)$ is called an {\em optimization space\/} if $X$ is a ground set and
$\CS$ is a set of trajectories on~$X$.

\subsection{Operational procedures.} \label{ss:operp}
If $(\Rb^+)^X$ denotes the set of all functions of the form $C:X\to\Rb^+$, we let
$\mathcal{C}(X)$ be a subset of $(\Rb^+)^X$, called the set of all (admissible)
{\em cost functions}, also representing distance, weight, time, etc. Given
$C\in\mathcal{C}(X)$, to each element $x\in X$ a nonnegative number $C(x)$
is assigned uniquely, which is called the {\em cost\/} of~$x$.

Since we are going to optimize (i.e., look for minima or maxima of) nonnegative functions
on the set of trajectories $\CS$, we denote by $\mbox{Ob}(\CS)=(\Rb^+)^\CS$ the
family of all such functions, called {\em objective functions}.

A map of the form $f:\mathcal{C}(X)\to\mbox{Ob}(\CS)$ is said to be an {\em
operational procedure\/} on the optimization space $(X,\CS)$. In other words, to each
cost function $C:X\to\Rb^+$ the operational procedure $f$ assignes in a unique way
the objective function of the form $f_C\equiv f(C):\CS\to\Rb^+$. If the cost function
$C$ is fixed (somehow), notation $f_C(S)$ will be employed in place of $f(C)(S)$,
where $S\in\CS$.

Of particular importance for the developments to follow are operational procedures
generated by {\sf A}-operations on $\Rb^+$, which are most often encountered
in practice.

Let $(X,\CS)$ be an optimization space, $C:X\to\Rb^+$---a cost function and
$\oplus$---an {\sf A}-operation on $\Rb^+$. Suppose the set function
$F_C:2^X\to\Rb^+$ is given by
  \begin{equation} \label{e:oppr}
F_C(S)=\bigoplus_{y\in S}C(y)\quad\mbox{for \,all}\quad S\in\dot2^X,\quad
\mbox{and}\quad F_C(\varnothing)=\varnothing,
  \end{equation}
where (cf.~\eqref{e:e3})
  \begin{equation} \label{e:sinn}
\bigoplus_{y\in S}C(y)=\bigoplus_{i=1}^nC(b(i))
  \end{equation}
for a bijection $b:\{1,\dots,n\}\to S$ with $n=|S|$ (since $\oplus$ satisfies axioms (A.1)
and (A.2), the right-hand side in \eqref{e:sinn} is independent of a bijection $b$ chosen).
By virtue of (A.1) and (A.2), $F_C$ is a finitely additive {\em measure\/} on $2^X$ with
respect to the {\sf A}-operation $\oplus$, that is,
  \begin{equation*}
\mbox{if $S_1,S_2\in2^X$ and $S_1\cap S_2=\varnothing$, then
$F_C(S_1\cup S_2)=F_C(S_1)\oplus F_C(S_2)$,}
  \end{equation*}
the term $F_C(\varnothing)=\varnothing$ being omitted (cf.~\eqref{e:e12}). The measure
$F_C$ will be called the {\em operational measure\/} corresponding to $C$ and~$\oplus$.

The operational procedure $f:\mathcal{C}(X)\to\mbox{Ob}(\CS)$ on the
optimization space $(X,\CS)$, {\em generated by the {\sf A}-op\-eration $\oplus$},
is given by
  \begin{equation} \label{e:fCS}
f_C(S)=F_C(S)=\bigoplus_{y\in S}C(y)\quad\mbox{for \,all}\quad S\in\CS\quad
\mbox{and}\quad C\in\mathcal{C}(X),
  \end{equation}
i.e., $f_C=F_C|_\CS$ is the restriction of measure $F_C$ to the set of trajectories~$\CS$
on~$X$.

\subsection{Remark.} \label{r:delta}
If $\oplus$ is an {\sf A}-operation of addition, then the measure $F_C$ can be expressed
as $F_C(S)=\oplus_{x\in X}C(x)\delta_x(S)$, where $\delta_x:2^X\to\{0,1\}$ is the Dirac
measure (or point mass) concentrated at $x\in X$, i.e., given $S\subset X$, one has
$\delta_x(S)=1$ if $x\in S$, and $\delta_x(S)=0$ if $x\notin S$.

\subsection{Examples of operational measures.} \label{ss:eom}
Here we follow the pattern of Section~\ref{ss:exul}. Let $C\in\mathcal{C}(X)$
and $S\in\dot2^X$.

(a) If $u\oplus v=(u^p+v^p)^{1/p}$ on $\Rb^+=[0,\infty)$ with $p>0$, then the
$p$-sum operational measure is given by
  \begin{equation*}
F_C(S)=\biggl(\sum_{y\in S}C(y)^p\biggr)^{1/p}.
  \end{equation*}

(b) If $u\oplus v=\max\{u,v\}$ on $\Rb^+=[0,\infty)$, then the max (or bottleneck)
operational measure is of the form
  \begin{equation*}
F_C(S)=\max_{y\in S}C(y)=\max\bigl\{C(y):y\in S\bigr\}.
  \end{equation*}

(c) If $u\oplus v=uv$ on $\Rb^+=(0,\infty)$, then the product operational measure
is given by
  \begin{equation*}
F_C(S)=\prod_{y\in S}C(y)\quad\mbox{if}\quad C(y)>0\quad\mbox{for \,all}\quad
y\in X.
  \end{equation*}

\subsection{Optimization problems.} \label{ss:OP}
The triple $(X,\CS,f)$, where $(X,\CS)$ is an optimization space and $f$ is an operational
procedure on $(X,\CS)$, determines a (Discrete) Optimization Problem (\DOP, for short),
which is formulated as follows: given a cost function $C:X\to\Rb^+$, {\em minimize or
maximize the objective function $f_C$ on $\CS$}, that is, look for solutions to the
following extremal problem:
  \begin{equation} \label{e:minmax}
f_C(S)\to\min\mbox{ (or \,$\max$)},\qquad S\in\CS.
  \end{equation}
The set of trajectories $\CS$ in \eqref{e:minmax} plays the role of the set of all
feasible (or admissible) solutions to the \DOP.

Throughout the paper we concentrate only on the minimization problem \eqref{e:minmax}
with the objective function of the form \eqref{e:fCS}, i.e., problem \eqref{e:DO},
where $\oplus$ is an {\sf A}-operation on $\Rb^+$. Since the formulation of the problem
\eqref{e:DO} depends on $X$, $\CS$, $\oplus$ and $C$, we will also refer to the problem
\eqref{e:DO} as \DOP$(X,\CS,\oplus,C)$.

Examples of concrete \DOP{}s including combinatorial \DOP{}s are the
well known traveling salesman problem, shortest path problem, assignment problem,
Steiner problem, machine sequencing problem, min-cut problem, and many other
problems on graphs, matroids, etc.\ (we refer to
\cite{JOGO12,Gal,Gord,Gordeev,Lawler,Libura,Ravi,Shier,Sotskov,Tarjan}).

\subsection{Optimal solutions.} \label{ss:opts}
Given an \DOP$(X,\CS,\oplus,C)$, we denote by
  \begin{equation} \label{e:S*}
\CS^*\equiv\CS^*_C=\{S^*\in\CS:\mbox{$f_C(S^*)\le f_C(S)$ for all $S\in\CS$}\}
  \end{equation}
the set of all {\em optimal solutions\/} to the \DOP{} \eqref{e:DO}. The collection $\CS^*$
and its elements depend on the cost function $C$; however, if $C$ is fixed (in a context),
then it will be convenient (and brief) not to show the dependence $\CS^*=\CS^*(C)$
explicitly. Note that, since the set of trajectories $\CS$ is finite, optimal solutions
$S^*\in\CS^*$ always exist, i.e., we have $\CS^*\ne\varnothing$.

The minimal value of $f_C$ on $\CS$, called the {\em optimal value\/} of the
\DOP$(X,\CS,\oplus,C)$, is determined uniquely and is independent of an optimal
solution $S^*\in\CS^*$, and it will be denoted by
  \begin{equation} \label{e:ov}
f_C(\CS^*)=\min_{S\in\CS}f_C(S)=\min_{\CS}f_C=f_C(S^*)\quad\mbox{for \,all}
\quad S^*\in\CS^*.
  \end{equation}

\subsection{Equivalent \DOP{}s.} \label{ss:eqvD}
Given an \DOP$(X,\CS,\oplus,C)$ of the form \eqref{e:DO} and a $\varphi$-function
$\varphi:[0,\infty)\to[0,\infty)$ (cf.\ Section~\ref{ss:equiv}), we let
$\widehat\oplus=E_\varphi(\oplus)$
(i.e., $u\widehat{\oplus}v=\varphi^{-1}\bigl(\varphi(u)\oplus\varphi(v)\bigr)$
for all $u,v\ge0$) and $\widehat C=\varphi^{-1}\circ C$ (i.e.,
$\widehat C(y)=\varphi^{-1}(C(y))$ for all $y\in X$). We are going to show that
the \DOP$(X,\CS,\oplus,C)$ is {\em equivalent\/} to the
\DOP$(X,\CS,\widehat\oplus,\widehat C)$ in the sense that their sets of optimal
solutions, denoted here by $\CS^*(\oplus)$ and $\CS^*(\widehat\oplus)$,
respectively, coincide. In fact, first we note that if $u_1,\dots,u_n\ge0$, then
$\widehat\oplus_{i=1}^nu_i=\varphi^{-1}\bigl(\oplus_{i=1}^n\varphi(u_i)\bigr)$,
$n\in\Nb$. It follows that, given $S\in\CS$, for the corresponding objective functions
$f_C^\oplus$ and $f_{\widehat C}^{\widehat\oplus}$ we have
  \begin{eqnarray*}
f_C^\oplus(S)&=&\bigoplus_{y\in S}C(y)
  =\bigoplus_{y\in S}\varphi\bigl(\varphi^{-1}(C(y))\bigr)
  =\bigoplus_{y\in S}\varphi\bigl(\widehat C(y)\bigr)=\\[3pt]
&=&\varphi\biggl(\varphi^{-1}\Bigl(\bigoplus_{y\in S}\varphi%
  \bigl(\widehat C(y)\bigr)\Bigr)\biggr)
  =\varphi\biggl(\,\,\widehat{\bigoplus_{y\in S}}\,\,\widehat C(y)\biggr)
  =\varphi\Bigl(f_{\widehat C}^{\widehat\oplus}(S)\Bigr).
  \end{eqnarray*}
Since $\varphi$ is strictly increasing, it follows from \eqref{e:S*} that
$\CS^*(\oplus)=\CS^*(\widehat\oplus)$.

\section{Upper stability intervals.} \label{s:usi}

\subsection{{}} \label{ss:sa}
In the Sensitivity Analysis one is interested in numerical characteristics of elements $x$
from the ground set $X$, which express the degree of invariance of an optimal solution
to the \DOP{} \eqref{e:DO} with respect to a change of the single cost $C(x)$.
The following two notions serve this purpose and are adopted in the literature
(\cite{DAN12}--\cite{Gordeev}, \cite{Libura}--\cite{Libura10}, \cite{Shier}--\cite{Tarjan}).
By the {\em upper tolerance\/} $u_{S^*}(x)$ ({\em lower tolerance $\ell_{S^*}(x)$})
of $x\in X$ with respect to $S^*\in\CS^*$ one means the {\em maximum increase\/}
({\em maximum decrease}, respectively) of the cost $C(x)$ only, so that the optimal
solution $S^*$ to the original \DOP{} \eqref{e:DO} remains an optimal solution to the
``perturbed'' \DOP{} \eqref{e:DO}, in which the costs $C(y)$ are unchanged if $y\ne x$
and the cost $C(x)$ of $x$ is increased (decreased, respectively) as compared to~$C(x)$.
These two notions will be studied in detail in the framework of general {\sf A}-operations
in this and the next sections. 

\smallbreak
Let the \DOP$(X,\CS,\oplus,C)$ of the form \eqref{e:DO} be given and $\CS^*$ be the
set of all its optimal solutions (cf.~\eqref{e:S*}).

\subsection{Perturbed objective functions.} \label{ss:pof}
Given $x\in X$, we perturb the cost function $C\in\mathcal{C}(X)$ at its value $C(x)$
by defining the {\em perturbed cost function} $C_{x,\gamma}:X\to\Rb^+$ with
$\gamma\in\Rb^+$ as follows: $C_{x,\gamma}(y)=C(y)$ if  $y\in X$ and $y\ne x$,
and $C_{x,\gamma}(x)=\gamma$. Since $C_{x,\gamma}\in\mathcal{C}(X)$, we let
$f(C_{x,\gamma})\equiv f_{C_{x,\gamma}}$ be the objective function \eqref{e:fCS}
corresponding to the cost function $C_{x,\gamma}$, called the {\em perturbed
objective function\/} (as compared to $f_C$), and so, it is of the form
  \begin{equation} \label{e:pof}
f(C_{x,\gamma})(S)=\bigoplus_{y\in S}C_{x,\gamma}(y)\quad\mbox{for \,all}
\quad S\in\CS.
  \end{equation}

In order to (properly) define upper and lower tolerances of $x\in X$ with respect to
an optimal solution $S^*\in\CS^*$ to problem \eqref{e:DO}, we ought to determine
(further) restrictions on $\gamma\in\Rb^+$, under which
  \begin{equation} \label{e:impl}
S^*\in\CS^*\quad\mbox{implies}\quad f(C_{x,\gamma})(S^*)
=\min_{S\in\CS}f(C_{x,\gamma})(S).
  \end{equation}
For this, let us express the perturbed objective function \eqref{e:pof} in terms of the
original cost function $C$, the initial objective function $f_C$ and the operational
measure~$F_C$.

In order to do it, it will be helpful to introduce two ad hoc subcollections of the set of
trajectories $\CS$ by
  \begin{equation} \label{e:Spmx}
\CS_{-x}=\{S\in\CS:x\notin S\}\quad\mbox{and}\quad\CS_x=\{S\in\CS:x\in S\},
\qquad x\in X; 
  \end{equation}
in other words, given $x\in X$ and $S\in\CS$, we have:
$x\notin S$ \,iff \,$S\in\CS_{-x}$, and $x\in S$ \,iff \,$S\in\CS_x$.
By virtue of \eqref{e:traj}, both collections $\CS_{-x}$ and $\CS_x$ are nonempty,
$\CS_{-x}\cup\CS_x=\CS$ and $\CS_{-x}\cap\CS_x=\varnothing$ for all $x\in X$.

Now, given $S\in\CS$, we have either $S\in\CS_{-x}$ or $S\in\CS_x$. If $S\in\CS_{-x}$
(or $x\notin S$), then $C_{x,\gamma}(y)=C(y)$ for all $y\in S$, and so,
\eqref{e:pof} and \eqref{e:fCS} imply
  \begin{equation*}
f(C_{x,\gamma})(S)=\bigoplus_{y\in S}C(y)=F_C(S)=f_C(S).
  \end{equation*}
If $S\in\CS_x$ (or $x\in S$), then $C_{x,\gamma}(y)=C(y)$ if $y\in S$ and $y\ne x$,
and $C_{x,\gamma}(x)=\gamma$, and so, \eqref{e:pof}, \eqref{e:oppr}, (A.1) and
(A.2) yield
  \begin{equation*}
f(C_{x,\gamma})(S)=\gamma\oplus\biggl(\bigoplus_{y\in S\setminus\{x\}}C(y)\biggr)
=\gamma\oplus F_C(S\setminus\{x\}),
  \end{equation*}
where the term $F_C(S\setminus\{x\})$ is omitted if $S=\{x\}$ (cf.~\eqref{e:e12}).

Thus, given $x\in X$ and $S\in\CS$, the perturbed objective function $f(C_{x,\gamma})$
is expressed as
  \begin{equation} \label{e:varpo}
f(C_{x,\gamma})(S)=\left\{
  \begin{array}{ccl}
f_C(S) & \mbox{if} & \mbox{$S\in\CS_{-x}$ (i.e., $x\notin S$),}\\[2pt]
\gamma\oplus F_C(S\setminus\{x\}) & \mbox{if} & 
    \mbox{$S\in\CS_x$ \,\,\,\,(i.e., $x\in S$).}
  \end{array}\right.
  \end{equation}
In particular, if $\gamma=C(x)$, then $C_{x,\gamma}(y)=C(y)$ for all $y\in X$, and so,
\eqref{e:varpo} implies
  \begin{equation}\label{e:fxs}
f_C(S)=C(x)\oplus F_C(S\setminus\{x\})\quad\mbox{if}\quad S\in\CS_x.
  \end{equation}

Formula \eqref{e:varpo} is valid for all {\sf A}-operations $\oplus$ on $\Rb^+$ and,
particularly, as it will be seen later, it works well for all strict operations.
However, for certain nonstrict {\sf A}-operations, such as $\max$, a different
(more subtle) form of formula \eqref{e:varpo} is needed.

\subsection{The perturbed objective function in the case $\oplus=\max$.}\label{ss:pmax}
Suppose that $\oplus=\max$ on $\Rb^+=[0,\infty)$. Let us transform the lower
part of formula \eqref{e:varpo} taking into account certain specific features of the
{\sf A}-operation~$\max$.

Only the case $S\in\CS_x$ with $S\ne\{x\}$ is to be considered. By virtue of \eqref{e:fxs}
and Example~\ref{ss:eom}(b), we have
  \begin{equation} \label{e:nfcs}
f_C(S)=\max\{C(x),F_C(S\setminus\{x\})\}\quad\mbox{with}\quad
F_C(S\setminus\{x\})=\max_{y\in S\setminus\{x\}}C(y),
  \end{equation}
and the second line of \eqref{e:varpo} can be rewritten as
  \begin{equation} \label{e:49}
f(C_{x,\gamma})(S)=\max\{\gamma,F_C(S\setminus\{x\})\}.
  \end{equation}

Assuming that $\gamma\ge C(x)$ and considering the two possibilities in \eqref{e:nfcs},
which are of the form
  \begin{equation} \label{e:mposs}
f_C(S)=C(x)\quad\mbox{or}\quad f_C(S)=F_C(S\setminus\{x\}),
  \end{equation}
we find from \eqref{e:49} that
  \begin{equation} \label{e:msx}
f(C_{x,\gamma})(S)=\max\{\gamma,f_C(S)\}=\gamma\oplus f_C(S).
  \end{equation}
In fact, if $f_C(S)=C(x)$, then
  \begin{equation*}
F_C(S\setminus\{x\})=\max_{y\in S\setminus\{x\}}C(y)\le\max_{y\in S}C(y)
=f_C(S)=C(x)\le\gamma,
  \end{equation*}
and so, by \eqref{e:49}, $f(C_{x,\gamma})(S)=\gamma$. Now, if
$f_C(S)=F_C(S\setminus\{x\})$, then once again \eqref{e:49} implies
$f(C_{x,\gamma})(S)=\max\{\gamma,f_C(S)\}$.

Assume that $0\le\gamma\le C(x)$ and consider the possibilities \eqref{e:mposs}.
If $f_C(S)=C(x)$, then, by \eqref{e:49}, we have
  \begin{equation} \label{e:forp}
f(C_{x,\gamma})(S)=\max\{f_C(S)+\gamma-C(x),F_C(S\setminus\{x\})\},
  \end{equation}
and if $f_C(S)=F_C(S\setminus\{x\})$, then
  \begin{equation} \label{e:gorn}
\gamma\le C(x)\le\max_{y\in S}C(y)=f_C(S)=F_C(S\setminus\{x\}),
  \end{equation}
and so, \eqref{e:49} implies $f(C_{x,\gamma})(S)=F_C(S\setminus\{x\})=f_C(S)$.
It follows that equality \eqref{e:forp} holds under both possibilities \eqref{e:mposs}.

Thus, given $x\in X$ and $S\in\CS$, taking into account \eqref{e:varpo}, \eqref{e:msx}
and \eqref{e:forp}, we have the following alternative expression for the perturbed
objective function in the case $\oplus=\max$:
  \begin{equation} \label{e:altma}
f(C_{x,\gamma})(S)=\left\{
  \begin{array}{ccl}
f_C(S) & \mbox{if} & S\in\CS_{-x},\\[2pt]
\max\{\gamma,f_C(S)\} & \mbox{if} & 
    \mbox{$S\in\CS_x$  and $\gamma\ge C(x)$,}\\[2pt]
\max\bigl\{f_C(S)\!+\!\gamma\!-\!C(x),F_C(S\!\setminus\!\{x\})\bigr\} &
    \mbox{if} & \mbox{$S\in\CS_x$ and $0\!\le\!\gamma\!\le\! C(x)$,}
  \end{array}\right.
  \end{equation}
where $f_C(S)=F_C(S)$ is as in Example~\ref{ss:eom}(b) and $F_C(S\setminus\{x\})$
is given in \eqref{e:nfcs}.

\subsection{Unrestricted upper tolerances.} \label{ss:uut}
In order to define the upper tolerance $u_{S^*}(x)$ of $x\in X$ with respect to an
$S^*\in\CS^*$ following the pattern exposed in Section~\ref{ss:sa}, we have to increase
the cost $C(x)$ to the value $\gamma\ge C(x)$ in such a way that the implication
\eqref{e:impl} holds, i.e., $S^*$ is also an optimal solution to the perturbed
\DOP$(X,\CS,\oplus,C_{x,\gamma})$, which is the problem \eqref{e:DO} with $C$
replaced by $C_{x,\gamma}$ (cf.\ also \eqref{e:pof}). It is to be noted that for
certain elements $x\in X$ implication \eqref{e:impl} holds automatically for all
$\gamma\ge C(x)$; for instance, it is intuitively clear that $x\notin S^*$ are such
elements, and so, the upper tolerance $u_{S^*}(x)$ for them may be thought of as
infinite (unrestricted). This assertion is made precise in the following

\subsection{Lemma.} \label{l:unut} {\em
If $S^*\in\CS^*$ and $x\in X\setminus S^*$, then for all $\gamma\ge C(x)$ we have
  \begin{equation} \label{e:fas}
f(C_{x,\gamma})(S^*)\le f(C_{x,\gamma})(S)\quad\mbox{for \,all}\quad S\in\CS.
  \end{equation}
}

\proof{Proof.}
Let us fix $\gamma\ge C(x)$ arbitrarily. Since $S^*\in\CS^*$, it follows from
\eqref{e:S*} that $S^*\in\CS$ and $f_C(S^*)\le f_C(S)$ for all $S\in\CS$.
Assumption $x\notin S^*$ is equivalent to $S^*\in\CS_{-x}$, and so,
\eqref{e:varpo} implies
  \begin{equation*}
f(C_{x,\gamma})(S^*)=f_C(S^*).
  \end{equation*}
Now, given $S\in\CS$, we have either $S\in\CS_{-x}$ or $S\in\CS_x$. If $S\in\CS_{-x}$,
then taking into account \eqref{e:varpo}, we get (for all $\gamma\in\Rb^+$)
  \begin{equation*}
f(C_{x,\gamma})(S^*)=f_C(S^*)\le f_C(S)=f(C_{x,\gamma})(S).
  \end{equation*}
If $S\in\CS_x$, then, by virtue of \eqref{e:fxs}, the monotonicity of $\oplus$ (cf.\
axiom (A.3)) and \eqref{e:varpo}, we find
  \begin{eqnarray*}
f(C_{x,\gamma})(S^*)&=&f_C(S^*)\le f_C(S)
  =C(x)\oplus F_C(S\setminus\{x\})\le\\[3pt]
&\le&\gamma\oplus F_C(S\setminus\{x\})=f(C_{x,\gamma})(S),
  \end{eqnarray*}
which was to be proved.
\Halmos\endproof

Lemma~\ref{l:unut} is valid for all {\sf A}-operations $\oplus$ and even those satisfying
only axioms (A.1)--(A.3). However, for lower tolerances, to be considered in
Section~\ref{s:lt}, the situation is more subtle (cf.\ Theorem~\ref{t:lts}).

\subsection{Upper stability intervals.} \label{ss:upsi}
Intuitively, for minimization problems \eqref{e:DO} the cost $C(x)$ of an $x\in S^*$
cannot be increased indefinitely so that \eqref{e:impl} holds. So, we are interested in
finding (restrictions on $\oplus$ and) the largest closed interval of costs
$[C(x),C_{S^*}^+(x)]$ with $C^+_{S^*}(x)\ge C(x)$, called the {\em upper stability
interval}, such that the implication \eqref{e:impl} is valid for all
$\gamma\in[C(x),C^+_{S^*}(x)]$.

Given $S^*\in\CS^*$ and $x\in S^*$, we set (cf.~\eqref{e:Spmx})
  \begin{equation} \label{e:C+}
C^+_{S^*}(x)=\max\Gamma_{x,S^*}(\CS_{-x}),
  \end{equation}
where, given a subcollection $\CS_1\subset\CS$ (usually, $\CS_1=\CS_{-x}$, $\CS_x$
or~$\CS$),
  \begin{equation} \label{e:Gam}
\Gamma_{x,S^*}(\CS_1)=\bigl\{\gamma\in\Rb^+:
\mbox{$f(C_{x,\gamma})(S^*)\le f(C_{x,\gamma})(S)$ for all $S\in\CS_1$}\bigr\}.
  \end{equation}
Note that, since $S^*\in\CS^*$ and $C_{x,\gamma}(y)=C(y)$ for all $y\in X$ if
$\gamma=C(x)$, it follows from \eqref{e:S*} that $C(x)\in\Gamma_{x,S^*}(\CS_1)$
(cf.\ also \eqref{e:fre} below).

In order to evaluate (and estimate) the quantity \eqref{e:C+} (see Theorem~\ref{t:Cpp}),
it will be convenient to apply the notation of the form \eqref{e:ov} for subcollections
$\CS_{-x}$ and $\CS_x$ from \eqref{e:Spmx}:
  \begin{equation} \label{e:CS}
f_C(\CS_{-x}^*)=\min_{S\in\CS_{-x}}f_C(S)\quad\mbox{and}\quad
f_C(\CS_x^*)=\min_{S\in\CS_x}f_C(S),\qquad x\in X
  \end{equation}
(to avoid ambiguities, we may explicitly set
$\CS_{-x}^*=(\CS_{-x})^*\ne(\CS^*)_{-x}$
and $\CS_x^*=(\CS_x)^*\ne(\CS^*)_x$).

\smallbreak
In the next theorem $\upsub$ and $\losub$ denote the upper and lower
subtractions for~$\oplus$.

\subsection{Theorem.} \label{t:Cpp} {\em
Given $S^*\in\CS^*$ and $x\in S^*$, we have\/{\rm:}
  \begin{itemize}
\item[{\rm(a)}] $C^+_{S^*}(x)=f_C(\CS_{-x}^*)\upsub F_C(S^*\setminus\{x\});$
\item[{\rm(b)}] $C(x)\le C_1^+(x)\le C^+_{S^*}(x)\le C_2^+(x)$ and
$C_1^+(x)\in\Gamma_{x,S^*}(\CS_{-x})$, where
  \begin{eqnarray}
C_1^+(x)&=&f_C(\CS_{-x}^*)\upsub\,[f_C(\CS^*)\upsub C(x)],
  \label{e:C1+}\\[3pt]
C_2^+(x)&=&f_C(\CS_{-x}^*)\upsub\,[f_C(\CS^*)\losub C(x)];
  \label{e:C2+}
  \end{eqnarray} 
\item[{\rm(c)}] if, in addition, $\oplus$ is strict on\/ $\Rb^+$ or $\oplus=\max$
  on $[0,\infty)$, then
  \begin{equation} \label{e:C1+2}
C^+_{S^*}(x)=C_1^+(x)=C_2^+(x)
=[C(x)\oplus f_C(\CS_{-x}^*)]\upsub f_C(\CS^*)
=\max\Gamma_{x,S^*}(\CS),
  \end{equation}
the implication\/ \eqref{e:impl} holds for all $\gamma\in[C(x),C^+_{S^*}(x)]$ and,
moreover, in the case $\oplus=\max$ we also have
$C^+_{S^*}(x)=f_C(\CS_{-x}^*)$.
  \end{itemize}
}
\medbreak
\proof{Proof.}
(a) Since $x\in S^*$ iff $S^*\in\CS_x$, it follows from \eqref{e:varpo} that, given
$S\in\CS_{-x}$ and $\gamma\in\Rb^+$,
  \begin{equation} \label{e:zig}
f(C_{x,\gamma})(S^*)=\gamma\oplus F_C(S^*\setminus\{x\})\quad\mbox{and}
\quad f(C_{x,\gamma})(S)=f_C(S).
  \end{equation}
Therefore, \eqref{e:Gam} with $\CS_1=\CS_{-x}$ and \eqref{e:CS} yield
  \begin{align}
\Gamma_{x,S^*}(\CS_{-x})&=\bigl\{\gamma\in\Rb^+:
  \mbox{$\gamma\oplus F_C(S^*\setminus\{x\})\le f_C(S)$
  for all $S\in\CS_{-x}$}\bigr\}=\nonumber \\[3pt]
&=\bigl\{\gamma\in\Rb^+:\gamma\oplus F_C(S^*\setminus\{x\})
  \le f_C(\CS_{-x}^*)\bigr\}. \label{e:Gaa}
  \end{align}
Since $S^*\in\CS^*$ and $x\in S^*$, \eqref{e:fxs} and \eqref{e:S*} imply
  \begin{equation} \label{e:fre}
C(x)\oplus F_C(S^*\setminus\{x\})=f_C(S^*)\le f_C(S)\quad\mbox{for \,all}\quad
S\in\CS_{-x},
  \end{equation}
and so,
  \begin{equation} \label{e:bz}
C(x)\oplus F_C(S^*\setminus\{x\})\le f_C(\CS_{-x}^*).
  \end{equation}
Hence, $C(x)\in\Gamma_{x,S^*}(\CS_{-x})$ and the pair
$(f_C(\CS_{-x}^*),F_C(S^*\setminus\{x\}))$ belongs to the domain $D(\upsub)$
of the upper subtraction $\upsub$ for~$\oplus$. Definitions \eqref{e:C+} and
\eqref{e:dups}, \eqref{e:Gaa} and \eqref{e:bz} imply inequality
$C^+_{S^*}(x)\ge C(x)$ and assertion~(a).

\smallbreak
(b) Taking into account the equality in \eqref{e:fre} and applying inequalities
\eqref{e:wvv} from Lemma~\ref{l:ineq} (with $w=F_C(S^*\setminus\{x\})$ and
$v=C(x)$), we get
  \begin{equation} \label{e:due}
f_C(S^*)\losub C(x)\le F_C(S^*\setminus\{x\})\le f_C(S^*)\upsub C(x),
  \end{equation}
where, by \eqref{e:ov}, $f_C(S^*)=f_C(\CS^*)$ is the optimal value of
problem~\eqref{e:DO}.

Now, we put $w=f_C(\CS_{-x}^*)$.

First, we set $u=f_C(S^*)\losub C(x)$ and $v=F_C(S^*\setminus\{x\})$. By
\eqref{e:due}, $u\le v$, and it follows from \eqref{e:bz} that $(w,v)\in D(\upsub)$,
and so, applying Lemma~\ref{l:scnde}(a), we get $(w,u)\in D(\upsub)$ and
$w\upsub v\le w\upsub u$. This inequality together with Theorem~\ref{t:Cpp}(a)
and \eqref{e:C2+} gives $C^+_{S^*}(x)\le C_2^+(x)$.

Second, we set $u=F_C(S^*\setminus\{x\})$ and $v=f_C(S^*)\upsub C(x)$. Then
\eqref{e:due} implies $u\le v$, and equality \eqref{e:ups} from Lemma~\ref{l:ineq}
and (A.2) give
  \begin{equation} \label{e:Cxv}
C(x)\oplus v=C(x)\oplus[f_C(S^*)\upsub C(x)]=f_C(S^*)\le f_C(\CS^*_{-x})=w,
  \end{equation}
i.e., $(w,v)\in D(\upsub)$. Applying Lemma~\ref{l:scnde}(a), we find
$w\upsub v\le w\upsub u$, and so, by \eqref{e:C1+} and Theorem~\ref{t:Cpp}(a),
$C_1^+(x)\le C^+_{S^*}(x)$.

In order to show that $C_1^+(x)\in\Gamma_{x,S^*}(\CS_{-x})$, we apply the notation
for $u$, $v$ and $w$ from the previous paragraph. By virtue of \eqref{e:C1+},
inequality $u\le v$, Lemma~\ref{l:scnde}(a) and \eqref{e:ups} from Lemma~\ref{l:ineq},
we get
  \begin{equation*}
C_1^+(x)\oplus F_C(S^*\setminus\{x\})=(w\upsub v)\oplus u\le
(w\upsub u)\oplus u=w=f_C(\CS_{-x}^*),
  \end{equation*}
and it remains to take into account \eqref{e:Gaa}.

The inequality $C(x)\le C_1^+(x)$ is a consequence of \eqref{e:Cxv},
\eqref{e:dups} and \eqref{e:C1+}.

\smallbreak
(c) First, we establish two auxiliary inequalities (under general conditions on~$\oplus$).
By virtue of Lemma~\ref{l:luineq}(a) and \eqref{e:C1+}, we find
  \begin{equation} \label{e:rir}
[C(x)\oplus f_C(\CS^*_{-x})]\losub f_C(\CS^*)\le
f_C(\CS_{-x}^*)\upsub\,[f_C(\CS^*)\upsub C(x)]=C_1^+(x),
  \end{equation}
and \eqref{e:C2+} and Lemma~\ref{l:luineq}(c) yield
  \begin{equation} \label{e:vav}
C_2^+(x)=f_C(\CS^*_{-x})\upsub\,[f_C(\CS^*)\losub C(x)]\le
[C(x)\oplus f_C(\CS_{-x}^*)]\upsub f_C(\CS^*).
  \end{equation}

1. Suppose that the {\sf A}-operation $\oplus$ is strict on $\Rb^+$. Then, by
\eqref{e:u=l} from Lemma~\ref{l:ineq}, $\upsub=\losub$ on $D(\upsub)$, and so,
the first three equalities in \eqref{e:C1+2} follow from Theorem~\ref{t:Cpp}(b),
\eqref{e:C1+}, \eqref{e:C2+}, \eqref{e:rir} and \eqref{e:vav}.

Now, given $\gamma\in[C(x),C^+_{S^*}(x)]$, let us show that \eqref{e:impl}
(or \eqref{e:fas}) holds. If $S\in\CS$, then either $S\in\CS_{-x}$ or $S\in\CS_x$.
Let $S\in\CS_{-x}$. Since $C(x)\le\gamma\le C_{S^*}^+(x)$, we have
$\gamma\in\Gamma_{x,S^*}(\CS_{-x})$, which implies the inequality in \eqref{e:fas}.
In more details, by \eqref{e:zig}, \eqref{e:Gaa}, \eqref{e:CS} and \eqref{e:varpo},
we have
  \begin{align}
f(C_{x,\gamma})(S^*)&=\gamma\oplus F_C(S^*\setminus\{x\})\le
  C^+_{S^*}(x)\oplus F_C(S^*\setminus\{x\})\le\nonumber\\[3pt]
&\le f_C(\CS_{-x}^*)\le f_C(S)=f(C_{x,\gamma})(S).\nonumber
  \end{align}
If $S\in\CS_x$ (and $\gamma\in\Rb^+$), then it follows from \eqref{e:fxs}
and \eqref{e:S*} that
  \begin{equation*}
C(x)\oplus F_C(S^*\setminus\{x\})=f_C(S^*)\le f_C(S)
=C(x)\oplus F_C(S\setminus\{x\}),
  \end{equation*}
and so, by the cancellation law \eqref{e:canl},
$F_C(S^*\setminus\{x\})\le F_C(S\setminus\{x\})$. Taking into account the
monotonicity of $\oplus$ and \eqref{e:varpo}, we obtain
  \begin{equation*}
f(C_{x,\gamma})(S^*)=\gamma\oplus F_C(S^*\setminus\{x\})\le
\gamma\oplus F_C(S\setminus\{x\})=f(C_{x,\gamma})(S).
  \end{equation*}
This proves also that $\Gamma_{x,S^*}(\CS_x)=\Rb^+$, whence
  \begin{equation*}
\Gamma_{x,S^*}(\CS)=\Gamma_{x,S^*}(\CS_{-x})\cap\Gamma_{x,S^*}(\CS_x)
=\Gamma_{x,S^*}(\CS_{-x}),
  \end{equation*}
and so, $C^+_{S^*}(x)=\max\Gamma_{x,S^*}(\CS)$, which is the fourth
equality in~\eqref{e:C1+2}.

2. Now assume that $\oplus=\max$ on $\Rb^+=[0,\infty)$. Since $x\in S^*$, by
\eqref{e:nfcs}, we have $C(x)\le f_C(S^*)=f_C(\CS^*)$, and so, taking into account
Example~\ref{ss:exul}(b), we get $f_C(\CS^*)\upsub C(x)=f_C(\CS^*)$, whereas
$f_C(\CS^*)\losub C(x)=f_C(\CS^*)$ if $C(x)<f_C(\CS^*)$, and
$f_C(\CS^*)\losub C(x)=0$ if $C(x)=f_C(\CS^*)$.
Since $S^*$ is an optimal solution to \eqref{e:DO}, we get
$f_C(\CS^*)\le f_C(\CS^*_{-x})$, and so, \eqref{e:C1+} implies
  \begin{equation*}
C_1^+(x)=f_C(\CS^*_{-x})\upsub f_C(\CS^*)=f_C(\CS^*_{-x}).
  \end{equation*}
By virtue of \eqref{e:C2+}, we find: if $C(x)<f_C(\CS^*)$, then
  \begin{equation*}
C_2^+(x)=f_C(\CS^*_{-x})\upsub f_C(\CS^*)=f_C(\CS^*_{-x}),
  \end{equation*}
and if $C(x)=f_C(\CS^*)$, then
  \begin{equation*}
C_2^+(x)=f_C(\CS^*_{-x})\upsub0=f_C(\CS^*_{-x}).
  \end{equation*}
Thus, $C_1^+(x)=C_2^+(x)=f_C(\CS^*_{-x})$ and, by Theorem~\ref{t:Cpp}(b),
$C^+_{S^*}(x)=f_C(\CS^*_{-x})$, which establishes the first two equalities
in~\eqref{e:C1+2}.

Since $C(x)\le f_C(\CS^*)\le f_C(\CS^*_{-x})$, we get (with $\oplus=\max$)
  \begin{equation*}
[C(x)\oplus f_C(\CS^*_{-x})]\upsub f_C(\CS^*)=\max\{C(x),f_C(\CS^*_{-x})\}
=f_C(\CS^*_{-x}),
  \end{equation*}
and the third equality in \eqref{e:C1+2} follows.

Now we prove \eqref{e:fas} for all $\gamma\in[C(x),C^+_{S^*}(x)]$. Since
$\gamma\le C^+_{S^*}(x)$, \eqref{e:fas} for $S\in\CS_{-x}$ follows from
the inclusion $\gamma\in\Gamma_{x,S^*}(\CS_{-x})$. In more details, given
$S\in\CS_{-x}$, \eqref{e:altma} implies
  \begin{align}
f(C_{x,\gamma})(S^*)&=\max\{\gamma,f_C(S^*)\}\le
  \max\{C^+_{S^*}(x),f_C(S^*)\}=\nonumber\\[3pt]
&=C^+_{S^*}(x)=f_C(\CS^*_{-x})\le f_C(S)=f(C_{x,\gamma})(S).\nonumber
  \end{align}
If $S\in\CS_x$ (and $\gamma\in\Rb^+$), then the monotonicity of $\max$ and
\eqref{e:altma} yield
  \begin{equation*}
f(C_{x,\gamma})(S^*)=\max\{\gamma,f_C(S^*)\}\le\max\{\gamma,f_C(S)\}
=f(C_{x,\gamma})(S).
  \end{equation*}
As in step~1 of item (c), this proves also that 
$C^+_{S^*}(x)=\max\Gamma_{x,S^*}(\CS)$.
\Halmos\endproof

\subsection{Remark.} \label{r:ios}
By Theorem~\ref{t:Cpp}(a) and definition \eqref{e:C+}, the value $C^+_{S^*}(x)$
depends on the optimal solution $S^*$ to problem \eqref{e:DO}. However, if $\oplus$
is strict or $\oplus=\max$, then, by Theorem~\ref{t:Cpp}(c), the value
$C^+_{S^*}(x)=C_1^+(x)$ is independent of optimal solutions $S^*$ to \eqref{e:DO}
such that $x\in S^*$ in the following sense: given $S_1^*,S_2^*\in\CS^*$, if
$x\in S_1^*\cap S_2^*$, then $C^+_{S^*_1}(x)=C^+_{S_2^*}(x)$ (in fact, these
two quantities are given by the same formula \eqref{e:C1+2}, which does not involve
neither $S_1^*$ nor~$S_2^*$).

\subsection{{}} Having the upper stability interval $[C(x),C^+_{S^*}(x)]$
(with $C(x)\le C^+_{S^*}(x)$) for $S^*\in\CS^*$ and $x\in S^*$, it may look
quite natural to define the {\em upper tolerance\/} $u_{S^*}(x)$ of $x\in S^*$
as a ``measure'' ($=$\,some generalized length) of the upper stability interval.
This can be done in many ways. For instance, if $\varphi:[0,\infty)\to[0,\infty)$
is a $\varphi$-function, then we may have set
$\mu([C(x),C^+_{S^*}(x)])=\varphi(C^+_{S^*}(x))-\varphi(C(x))$.
However, this is irrelevant for our purposes, because the upper stability interval has
been generated via the {\sf A}-operation $\oplus$ or, more precisely (cf.\
Theorem~\ref{t:Cpp}(a)), by the upper subtraction $\upsub$ for~$\oplus$.
Having this in mind, as well as the translation invariance of $\upsub$ and $\losub$
(cf.\ Lemma~\ref{l:trinv}), we adopt the following definition.

\subsection{Definition.} \label{d:utol}
Assume that the {\sf A}-operation $\oplus$ is strict on $\Rb^+$. Given $S^*\in\CS^*$
and $x\in S^*$, the {\em upper tolerance\/} of $x$ is defined by
  \begin{equation} \label{e:u}
u_{S^*}(x)=C^+_{S^*}(x)\upsub C(x)\in\Rb^+.
  \end{equation}

\subsection{Theorem.} \label{t:upto} {\em
If\/ $\oplus$ is a strict\/ {\sf A}-operation on $\Rb^+$, $S^*\in\CS^*$ and $x\in S^*$,
then the value $u_{S^*}(x)$ is well-defined,
  \begin{equation} \label{e:ucal}
u_{S^*}(x)=f_C(\CS^*_{-x})\upsub f_C(\CS^*)\ge{\pmb e},
  \end{equation}
and $u_{S^*}(x)={\pmb e}$ iff\/ $C^+_{S^*}(x)=C(x)$ iff\/ 
$f_C(\CS^*_{-x})=f_C(\CS^*)$, where ${\pmb e}\in\Rb^+$ is the neutral element
with respect to~$\oplus$.
}
\medbreak
\proof{Proof.}
If $\oplus$ is a generalized addition (i.e., satisfies (A.1)--(A.3) and (A.5)), then
inequality $C(x)\le C^+_{S^*}(x)$ from Theorem~\ref{t:Cpp}(b),\,(c) and
Lemma~\ref{l:comp}(b) imply $(C(x),C^+_{S^*}(x))\in D(\upsub)$, and if $\oplus$
is a generalized multiplication (i.e., (A.1)--(A.4) and (A.6) are satisfied), then the same
inclusion is a consequence of Lemma~\ref{l:comp}(c). It follows from \eqref{e:dups}
that the quantity \eqref{e:u} is well-defined.

Since ${\pmb e}\oplus C(x)=C(x)\le C^+_{S^*}(x)$, definition \eqref{e:dups} implies
  \begin{equation*}
{\pmb e}\le C^+_{S^*}(x)\upsub C(x)=u_{S^*}(x).
  \end{equation*}

If $u_{S^*}(x)={\pmb e}$, then, by virtue of \eqref{e:ups}, we get
  \begin{equation*}
C^+_{S^*}(x)=(C^+_{S^*}(x)\upsub C(x))\oplus C(x)=u_{S^*}(x)\oplus C(x)
={\pmb e}\oplus C(x)=C(x).
  \end{equation*}
(Note that all assertions above do not rely on the strictness of $\oplus$.) Now,
if equality $C^+_{S^*}(x)=C(x)$ holds, then we claim that $u_{S^*}(x)={\pmb e}$,
for, otherwise, if $u_{S^*}(x)>{\pmb e}$, then \eqref{e:ups} and the strictness
(A.3${}_{\mbox{\footnotesize\rm s}}$) of $\oplus$ imply
  \begin{equation*}
C^+_{S^*}(x)=(C^+_{S^*}(x)\upsub C(x))\oplus C(x)=u_{S^*}(x)\oplus C(x)
>{\pmb e}\oplus C(x)=C(x),
  \end{equation*}
which contradicts the assumption.

Finally, let us establish the equality in \eqref{e:ucal}. Setting
  \begin{equation*}
u=C(x),\quad v=f_C(\CS^*)\quad\mbox{and}\quad w=f_C(\CS^*_{-x}),
  \end{equation*}
we find from \eqref{e:u} and \eqref{e:C1+2} that
  \begin{equation*}
u_{S^*}(x)=C^+_{S^*}(x)\upsub C(x)=\bigl([u\oplus w]\upsub v\bigr)\upsub u.
  \end{equation*}
By \eqref{e:ups}, $w=(w\upsub v)\oplus v$, and so, \eqref{e:u=l} and
\eqref{e:wvv} yield
  \begin{equation*}
[u\oplus w]\upsub v=[u\oplus(w\upsub v)\oplus v]\upsub v=u\oplus(w\upsub v)
  \end{equation*}
and
  \begin{equation*}
\bigl([u\oplus w]\upsub v\bigr)\upsub u=\bigl(u\oplus(w\upsub v)\bigr)\upsub u
=w\upsub v,
  \end{equation*}
and it remains to take into account that
$w\upsub v=f_C(\CS^*_{-x})\upsub f_C(\CS^*)$.
\Halmos\endproof

\subsection{Examples of upper stability intervals and upper tolerances.} \label{ss:ust}
In accordance with Examples~\ref{ss:exul}(a) and (c), \eqref{e:C1+2} and
\eqref{e:ucal}, given $S^*\in\CS^*$ and $x\in S^*$, we have:

(a) if $u\oplus v=(u^p+v^p)^{1/p}$ on $\Rb^+=[0,\infty)$ with $p>0$, then
  \begin{align}
C^+_{S^*}(x)&=\bigl(C(x)^p+f_C(\CS^*_{-x})^p-f_C(\CS^*)^p\bigr)^{1/p},
  \nonumber\\[2pt]
u_{S^*}(x)&=\bigl(f_C(\CS^*_{-x})^p-f_C(\CS^*)^p\bigr)^{1/p}.\nonumber
  \end{align}
In particular, if $p=1$, then for the \DOP{} $f_C(S)=\sum_{y\in S}C(y)\to\min$ with
$S\in\CS$ and an element $x\in S^*$ from its optimal solution $S^*\in\CS^*$ we find
  \begin{equation*}
u_{S^*}(x)=f_C(\CS^*_{-x})-f_C(\CS^*)
=\min_{\mbox{$\scriptstyle S\in\CS$\,\footnotesize\rm with\,$x\!\notin\! S$}}f_C(S)-
\min_{S\in\CS}f_C(S),
  \end{equation*}
which gives a formula due to Libura \cite{Libura} (cf.\ also \cite{JOGO12}, \cite{GJM},
\cite{Shier}--\cite{Tarjan});

\smallbreak
(b) if $u\oplus v=uv$ on $\Rb^+$ and $C(y)>0$ for all $y\in X$, then
  \begin{equation*}
C^+_{S^*}(x)=\displaystyle\frac{C(x)\!\cdot\! f_C(\CS^*_{-x})}{f_C(\CS^*)}
\,\,\quad\mbox{and}\quad\,\,
u_{S^*}(x)=\frac{f_C(\CS^*_{-x})}{f_C(\CS^*)}.
  \end{equation*}

\section{Lower stability intervals.} \label{s:lt}

Let the \DOP$(X,\CS,\oplus,C)$ of the form \eqref{e:DO} be given and $\CS^*$ be
the set of its optimal solutions from~\eqref{e:S*}.

In order to define the lower tolerance $\ell_{S^*}(x)$ of an $x\in X$ with respect to
an $S^*\in\CS^*$ (cf.\ Section~\ref{ss:sa}), we have to decrease the cost $C(x)$
to the value $\gamma\le C(x)$ so that $S^*$ is also an optimal solution to the
perturbed \DOP$(X,\CS,\oplus,C_{x,\gamma})$, i.e., the implication \eqref{e:impl} holds.

\subsection{Lower stability intervals.} \label{ss:losi}
One may (intuitively) feel that for minimization problems \eqref{e:DO} the cost of an
element $x\in X\setminus S^*$ cannot be decreased ``unboundedly'' in such a way
that the implication \eqref{e:impl} is valid. It is our aim now to obtain (restrictions on the
{\sf A}-operation $\oplus$ and) the largest closed interval of costs
$[C^-_{S^*}(x),C(x)]$ with $C^-_{S^*}(x)\le C(x)$, termed the {\em lower stability
interval}, such that \eqref{e:impl} holds for all $\gamma\in[C^-_{S^*}(x),C(x)]$.

Given $S^*\in\CS^*$ and $x\in X\setminus S^*$, we set (cf.\ \eqref{e:Spmx}
and \eqref{e:Gam})
  \begin{equation} \label{e:Cm}
C^-_{S^*}(x)=\min\Gamma_{x,S^*}(\CS_x).
  \end{equation}

In order to estimate and/or evaluate $C^-_{S^*}(x)$ (see Theorem~\ref{t:Cmm}),
we apply notation \eqref{e:CS} and note that there exists an $S_x\in\CS_x$ (i.e.,
$S_x\in\CS$ and $x\in S_x$) such that
  \begin{equation} \label{e:Sx}
F_C(S_x\setminus\{x\})=\min_{S\in\CS_x}F_C(S\setminus\{x\})
  \end{equation}
(actually, the value at the right in \eqref{e:Sx} and, hence, the quantity at the left
in \eqref{e:Sx}, are independent of the set $S_x$). Moreover, we have
  \begin{equation} \label{e:Sxeq}
f_C(S_x)=f_C(\CS^*_x).
  \end{equation}  
In fact, since $S_x\in\CS_x$, \eqref{e:CS} implies $f_C(S_x)\ge f_C(\CS^*_x)$;
on the other hand, given $S\in\CS_x$, since $x\in S$ and
$F_C(S_x\setminus\{x\})\le F_C(S\setminus\{x\})$, we find, by \eqref{e:fxs} and (A.3),
  \begin{align}
f_C(S_x)&=C(x)\oplus F_C(S_x\setminus\{x\})\le \label{e:5.4}\\[3pt]
&\le C(x)\oplus F_C(S\setminus\{x\})=f_C(S),\nonumber
  \end{align}
and so, again by \eqref{e:CS}, $f_C(S_x)\le f_C(\CS^*_x)$, which yields
equality~\eqref{e:Sxeq}.

In the next theorem $\upsub$ and $\losub$ denote (as usual) the upper and lower
subtractions for~$\oplus$, respectively.

\subsection{Theorem.} \label{t:Cmm} {\em
Given $S^*\in\CS^*$ and $x\in X\setminus S^*$, we have\/{\rm:}
  \begin{itemize}
\item[{\rm(a)}] $C^-_{S^*}(x)=f_C(\CS^*)\losub F_C(S_x\setminus\{x\})%
=\min\Gamma_{x,S^*}(\CS)$ {\rm(}cf.\ \eqref{e:Sx}{\rm)}, and the implication\/
\eqref{e:impl} holds for all $\gamma\in[C^-_{S^*}(x),C(x)];$
\item[{\rm(b)}] $C_1^-(x)\le C^-_{S^*}(x)\le C_2^-(x)\le C(x)$ and
$C_2^-(x)\in\Gamma_{x,S^*}(\CS)$, where
  \begin{align}
C_1^-(x)&=f_C(\CS^*)\losub\,[f_C(\CS^*_x)\upsub C(x)],\label{e:C1m}\\[3pt]
C_2^-(x)&=f_C(\CS^*)\losub\,[f_C(\CS^*_x)\losub C(x)];\label{e:C2m}
  \end{align}
\item[{\rm(c)}] if the {\sf A}-operation $\oplus$ is strict on $\Rb^+$, then
  \begin{equation} \label{e:Cmeq}
C^-_{S^*}(x)=C_1^-(x)=C_2^-(x)=[C(x)\oplus f_C(\CS^*)]\losub f_C(\CS^*_x);
  \end{equation}
\item[{\rm(d)}] if\/ $\oplus=\max$ on\/ $\Rb^+=[0,\infty)$, then $C_1^-(x)=0$,
  \begin{equation} \label{e:Ox}
C^-_{S^*}(x)=\left\{
  \begin{array}{ccl} 
f_C(\CS^*) & \mbox{if\/} & F_C(S_x\setminus\{x\})<f_C(\CS^*),\\[2pt]
0               & \mbox{if\/} & F_C(S_x\setminus\{x\})\ge f_C(\CS^*),
  \end{array}\right.
  \end{equation}
and
  \begin{equation} \label{e:cro}
C^-_2(x)=\left\{
  \begin{array}{ccl} 
f_C(\CS^*) & \mbox{if\/} & C(x)=f_C(\CS^*_x),\\[2pt]
0               & \mbox{if\/} & C(x)<f_C(\CS^*_x).
  \end{array}\right.
  \end{equation}
  \end{itemize}
}
\medbreak
\proof{Proof.}
(a) Since $S^*\in\CS_{-x}$, given $S\in\CS_x$ and $\gamma\in\Rb^+$, by
\eqref{e:varpo}, we find
  \begin{equation*}
f(C_{x,\gamma})(S^*)=f_C(S^*)=f_C(\CS^*)\quad\mbox{and}\quad
f(C_{x,\gamma})(S)=\gamma\oplus F_C(S\setminus\{x\}).
  \end{equation*}
It follows from \eqref{e:Gam} with $\CS_1=\CS_x$ and \eqref{e:Sx} that
  \begin{align}
\Gamma_{x,S^*}(\CS_x)&=\bigl\{\gamma\in\Rb^+:
  \mbox{$f_C(S^*)\le\gamma\oplus F_C(S\setminus\{x\})$ for all $S\in\CS_x$}\bigr\}
  =\nonumber\\[3pt]
&=\bigl\{\gamma\in\Rb^+:f_C(\CS^*)\le
\gamma\oplus F_C(S_x\setminus\{x\})\bigr\}. \label{e:hyp}
  \end{align}
Since $S^*\in\CS^*$, $f_C(S^*)\le f_C(S_x)$, and so, \eqref{e:5.4} implies
$C(x)\in\Gamma_{x,S^*}(\CS_x)$. Now, \eqref{e:Cm}, \eqref{e:hyp} and
definition \eqref{e:dlos} of the lower subtraction $\losub$ for~$\oplus$ yield inequality
$C^-_{S^*}(x)\le C(x)$ and the first equality in~(a).

Given $\gamma\in[C^-_{S^*}(x),C(x)]$, let us show that \eqref{e:impl} (or
\eqref{e:fas}) holds. For this, suppose $S\in\CS$. If $S\in\CS_{-x}$, then, taking into
account that $S^*\in\CS_{-x}\cap\CS^*$ and \eqref{e:varpo}, we have
(even for all $\gamma\in\Rb^+$)
  \begin{equation*}
f(C_{x,\gamma})(S^*)=f_C(S^*)\le f_C(S)=f(C_{x,\gamma})(S).
  \end{equation*}
This proves also that $\Gamma_{x,S^*}(\CS_{-x})=\Rb^+$, whence
  \begin{equation*}
\Gamma_{x,S^*}(\CS)=\Gamma_{x,S^*}(\CS_{-x})\cap\Gamma_{x,S^*}(\CS_x)
=\Gamma_{x,S^*}(\CS_x),
  \end{equation*}
and so, the second equality in (a) is a consequence of definition~\eqref{e:Cm}.

If $S\in\CS_x$, then inequalities $C^-_{S^*}(x)\le\gamma\le C(x)$ imply
$\gamma\in\Gamma_{x,S^*}(\CS_x)$, which gives the inequality in \eqref{e:fas}.
More directly, by virtue of \eqref{e:hyp}, \eqref{e:Sx} and \eqref{e:varpo},
  \begin{align}
f(C_{x,\gamma})(S^*)&=f_C(S^*)\le C^-_{S^*}(x)\oplus F_C(S_x\setminus\{x\})
  \le\gamma\oplus F_C(S_x\setminus\{x\})=\nonumber\\[2pt]
&=\gamma\oplus\Bigl(\min_{S'\in\CS_x}F_C(S'\setminus\{x\})\Bigr)
  \le\gamma\oplus F_C(S\setminus\{x\})=f(C_{x,\gamma})(S).\nonumber
  \end{align}

(b) Taking into account \eqref{e:Sxeq} and \eqref{e:5.4} and applying inequalities
\eqref{e:wvv} from Lemma~\ref{l:ineq} (with $w=F_C(S_x\setminus\{x\})$ and
$v=C(x)$), we find
  \begin{equation} \label{e:fine}
f_C(\CS^*_x)\losub C(x)\le F_C(S_x\setminus\{x\})\le f_C(\CS^*_x)\upsub C(x).
  \end{equation}

Now, let us put $w=f_C(\CS^*)$.

First, we set $u=f_C(\CS^*_x)\losub C(x)$ and $v=F_C(S_x\setminus\{x\})$.
By virtue of \eqref{e:fine}, $u\le v$, and so, Lemma~\ref{l:scnde}(b) implies
$w\losub v\le w\losub u$. This inequality, Theorem~\ref{t:Cmm}(a) and \eqref{e:C2m}
give $C^-_{S^*}(x)\le C_2^-(x)$. Since $C^-_{S^*}(x)$ is the minimal element of
the set $\Gamma_{x,S^*}(\CS)$, the last inequality yields
$C_2^-(x)\in\Gamma_{x,S^*}(\CS)$. It is also worth-while to verify this inclusion
directly: in fact, by virtue of \eqref{e:los}, (A.3), \eqref{e:C2m} and \eqref{e:hyp},
we have
  \begin{equation*}
f_C(\CS^*)=w\le(w\losub v)\oplus v\le(w\losub u)\oplus v
=C_2^-(x)\oplus F_C(S_x\setminus\{x\}).
  \end{equation*}
To see that $C_2^-(x)\le C(x)$, we note that, by \eqref{e:los},
  \begin{equation*}
f_C(\CS^*_x)\le C(x)\oplus[f_C(\CS^*_x)\losub C(x)],
  \end{equation*}
and so, the desired inequality follows from \eqref{e:dlos} and \eqref{e:C2m}.

Second, we set $u=F_C(S_x\setminus\{x\})$ and $v=f_C(\CS^*_x)\upsub C(x)$.
Then \eqref{e:fine} implies $u\le v$, and so, by Lemma~\ref{l:scnde}(b),
$w\losub v\le w\losub u$, and it follows from \eqref{e:C1m} and
Theorem~\ref{t:Cmm}(a) that $C_1^-(x)\le C^-_{S^*}(x)$.

\smallbreak
(c) By the strictness of $\oplus$ and \eqref{e:u=l}, $\upsub=\losub$ on $D(\upsub)$,
and so, \eqref{e:C1m}, \eqref{e:C2m} and inequalities in item (b) imply
$C_1^-(x)=C^-_{S^*}(x)=C_2^-(x)$, which proves the first two equalities
in~\eqref{e:Cmeq}.

In order to prove the third equality in \eqref{e:Cmeq}, we set
  \begin{equation} \label{e:doog}
C_0^-(x)=[C(x)\oplus f_C(\CS^*)]\losub f_C(\CS^*_x)
  \end{equation}
and note that, by virtue of Lemma~\ref{l:luineq}(b) and \eqref{e:C1m}, we have
(the following inequality, which is independent of the strictness of~$\oplus$)
  \begin{equation}\label{e:C0}
C_0^-(x)\le f_C(\CS^*)\losub\,[f_C(\CS^*_x)\upsub C(x)]=C_1^-(x).
  \end{equation}
Had we shown that $C_0^-(x)\in\Gamma_{x,S^*}(\CS_x)$, then, by \eqref{e:Cm},
we would have $C_1^-(x)=C^-_{S^*}(x)\le C_0^-(x)$, which implies the third
equality in~\eqref{e:Cmeq}. Setting $u=C(x)$, $v=F_C(S_x\setminus\{x\})$ and
$w=f_C(\CS^*)$ and noting that, by virtue of \eqref{e:Sxeq} and \eqref{e:5.4},
$f_C(\CS^*_x)=u\oplus v$, we find
  \begin{equation} \label{e:opop}
C_0^-(x)\oplus F_C(S_x\setminus\{x\})=\bigl[(u\oplus w)\losub(u\oplus v)\bigr]\oplus v.
  \end{equation}
If $u_1=(u\oplus w)\losub(u\oplus v)$, then inequality \eqref{e:los} implies
  \begin{equation*}
u\oplus w\le[(u\oplus w)\losub(u\oplus v)]\oplus(u\oplus v)=u_1\oplus u\oplus v,
  \end{equation*}
and so, cancelling by $u$ (by the strictness of $\oplus$ and \eqref{e:canl}), we get
$w\le u_1\oplus v$. It follows from definition \eqref{e:dlos} of $\losub$ that
$w\losub v\le u_1$. Hence, \eqref{e:los}, (A.3) and \eqref{e:opop} yield
  \begin{equation*}
f_C(\CS^*)=w\le(w\losub v)\oplus v\le u_1\oplus v
=C_0^-(x)\oplus F_C(S_x\setminus\{x\}),
  \end{equation*}
which, by virtue of \eqref{e:hyp}, gives $C_0^-(x)\in\Gamma_{x,S^*}(\CS_x)$.

\smallbreak
(d) Suppose $\oplus=\max$ on $\Rb^+=[0,\infty)$. Let us evaluate the values of
$C^-_{S^*}(x)$, $C_1^-(x)$ and $C_2^-(x)$ from Theorem~\ref{t:Cmm}(b) and the
value of $C_0^-(x)$ from \eqref{e:doog}. Note that, by virtue of item (b) and
\eqref{e:C0}, we have
  \begin{equation} \label{e:sh}
C_0^-(x)\le C_1^-(x)\le C^-_{S^*}(x)\le C_2^-(x)\le C(x).
  \end{equation}

First, equality \eqref{e:Ox} is a straightforward consequence of the first equality
in Theorem~\ref{t:Cmm}(a) and Example~\ref{ss:exul}(b).

Second, in order to evaluate $C_1^-(x)$, we note that $C(x)\le f_C(\CS^*_x)$;
in fact, if $S\in\CS_x$, then $x\in S$, and so, $C(x)\le\max_{y\in S}C(y)=f_C(S)$,
which,  by \eqref{e:CS}, implies the desired inequality. It follows from
Example~\ref{ss:exul}(b) that
  \begin{equation*}
f_C(\CS^*_x)\upsub C(x)=f_C(\CS^*_x),
  \end{equation*}
and, since $f_C(\CS^*)\le f_C(\CS^*_x)$, \eqref{e:C1m} and Example~\ref{ss:exul}(b)
give
  \begin{equation*}
C_1^-(x)=f_C(\CS^*)\losub\,[f_C(\CS^*_x)\upsub C(x)]=
f_C(\CS^*)\losub f_C(\CS^*_x)=0.
  \end{equation*}

Third, by the previous step and \eqref{e:sh}, we get $C_0^-(x)=0$. Also, this can be
seen directly as follows: since $C(x)\le f_C(\CS^*_x)$ and $f_C(\CS^*)\le f_C(\CS^*_x)$,
we find
  \begin{equation*}
C(x)\oplus f_C(\CS^*)=\max\{C(x),f_C(\CS^*)\}\le f_C(\CS^*_x),
  \end{equation*}
and so, \eqref{e:doog} and Example~\ref{ss:exul}(b) imply $C_0^-(x)=0$.

Fourth, since $C(x)\le f_C(\CS^*_x)$, it follows from Example~\ref{ss:exul}(b) that
  \begin{equation*}
f_C(\CS^*_x)\losub C(x)=\left\{
  \begin{array}{clc}
0 & \mbox{if} & C(x)=f_C(\CS^*_x),\\[2pt]
f_C(\CS^*_x) & \mbox{if} & C(x)<f_C(\CS^*_x),
  \end{array}\right.
  \end{equation*}
and so, \eqref{e:C2m} yields
  \begin{equation*}
C_2^-(x)=\left\{
  \begin{array}{clc}
f_C(\CS^*)\losub0 & \mbox{if} & C(x)=f_C(\CS^*_x),\\[2pt]
f_C(\CS^*)\losub f_C(\CS^*_x) & \mbox{if} & C(x)<f_C(\CS^*_x),
  \end{array}\right.
  \end{equation*}
whence \eqref{e:cro} follows if we take into account that $f_C(\CS^*)\le f_C(\CS^*_x)$.
\Halmos\endproof

\subsection{Remark.} \label{r:C-x}
By (the first equality in) Theorem~\ref{t:Cmm}(a), the value $C^-_{S^*}(x)$ with
$x\in X\setminus S^*$ does not depend on the optimal solution $S^*$ to problem
\eqref{e:DO} in the following sense: if $S_1^*,S_2^*\in\CS^*$ and
$x\in(X\setminus S_1^*)\cap(X\setminus S_2^*)$, then
$C^-_{S_1^*}(x)=C^-_{S_2^*}(x)$.

\smallbreak
In our next result we treat the case of ``unrestricted'' lower tolerances: for certain
elements $x$ from $X$ (e.g., $x\in S^*$) the implication \eqref{e:impl} always holds
for all costs $\gamma\le C(x)$. In particular, this clarifies definition \eqref{e:Cm}.
However, in contrast with Lemma~\ref{l:unut}, we will have to assume that the
{\sf A}-operation $\oplus$ is strict or $\oplus=\max$.

\subsection{Theorem.} \label{t:lts} {\em
Given $S^*\in\CS^*$ and $x\in S^*$, if one of the following two conditions\/ {\rm(a)}
or\/ {\rm(b)} holds\/{\rm:}
  \begin{itemize}
\item[{\rm(a)}] the\/ {\sf A}-operation $\oplus$ is strict on\/ $\Rb^+$,  or
\item[{\rm(b)}]
$\oplus=\max$ on\/ $[0,\infty)$ and either
  \begin{itemize}
\item[{\rm(i)}] $C(x)\!<\!f_C(\CS^*)$, or
\item[{\rm(ii)}] $C(x)\!=\!f_C(\CS^*)$ and $\CS^*\!=\!\{S^*\}$, or
\item[{\rm(iii)}] $C(x)\!=\!f_C(\CS^*)\!\le\! F_C(S_x\!\setminus\!\{x\})$, or
\item[{\rm(iv)}] $f_C(\CS^*)\!\le\! F_C(S_x\!\setminus\!\{x\})$,
  \end{itemize}
  \end{itemize}
then $f(C_{x,\gamma})(S^*)\le f(C_{x,\gamma})(S)$ for all $S\in\CS$ and
$\gamma\in\Rb^+$ with $\gamma\le C(x)$.
}
\medbreak
\proof{Proof.}
(a) Let $\oplus$ be strict and $\gamma\in\Rb^+$, $\gamma\le C(x)$, be arbitrarily
fixed. By \eqref{e:S*}, $S^*\in\CS$ and $f_C(S^*)\le f_C(S)$ for all $S\in\CS$.
Since $x\in S^*$ iff $S^*\in\CS_x$, \eqref{e:varpo} implies
$f(C_{x,\gamma})(S^*)=\gamma\oplus F_C(S^*\setminus\{x\})$. Given $S\in\CS$,
we have either $S\in\CS_{-x}$ or $S\in\CS_x$. If $S\in\CS_{-x}$, then the monotonicity
(A.3) of $\oplus$, \eqref{e:fxs} and \eqref{e:varpo} yield
  \begin{align*}
f(C_{x,\gamma})(S^*)&=\gamma\oplus F_C(S^*\setminus\{x\})\le
  C(x)\oplus F_C(S^*\setminus\{x\})=\\[3pt]
&=f_C(S^*)\le f_C(S)=f(C_{x,\gamma})(S).
  \end{align*}
Now, let $S\in\CS_x$, i.e., $x\in S$. Since, by \eqref{e:fxs} and \eqref{e:S*},
  \begin{equation*}
C(x)\oplus F_C(S^*\setminus\{x\})=f_C(S^*)\le f_C(S)
=C(x)\oplus F_C(S\setminus\{x\}),
  \end{equation*}
the strictness of $\oplus$ and the cancellation law \eqref{e:canl} imply
$F_C(S^*\setminus\{x\})\le F_C(S\setminus\{x\})$, and so, by (A.3) and
\eqref{e:varpo}, we get (even for all $\gamma\in\Rb^+$)
  \begin{equation*}
f(C_{x,\gamma})(S^*)=\gamma\oplus F_C(S^*\setminus\{x\})\le
\gamma\oplus F_C(S\setminus\{x\})=f(C_{x,\gamma})(S).
  \end{equation*}

(b) Suppose $\oplus=\max$ on $[0,\infty)$, and let us fix $0\le\gamma\le C(x)$.
Since $S^*\in\CS^*$, $f_C(S^*)\le f_C(S)$ for all $S\in\CS$. Because $x\in S^*$ iff
$S^*\in\CS_x$, \eqref{e:altma} gives
  \begin{equation} \label{e:mxy}
f(C_{x,\gamma})(S^*)=\max\{f_C(S^*)+\gamma-C(x),F_C(S^*\setminus\{x\})\}.
  \end{equation}
Let $S\in\CS$. If $S\in\CS_{-x}$, then, by \eqref{e:altma},
$f(C_{x,\gamma})(S)=f_C(S)$. The inequality $f_C(S^*)\le f_C(S)$ can be rewritten as
  \begin{equation*}
f_C(S^*)=\max\{C(x),F_C(S^*\setminus\{x\})\}\le f_C(S),
  \end{equation*}
which implies $F_C(S^*\setminus\{x\})\le f_C(S)$. Inequality $\gamma\le C(x)$
implies (cf.\ \eqref{e:mxy})
  \begin{equation} \label{e:cot}
f_C(S^*)+\gamma-C(x)\le f_C(S)+\gamma-C(x)\le f_C(S),
  \end{equation}
and so, by \eqref{e:mxy},
  \begin{equation*}
f(C_{x,\gamma})(S^*)\le f_C(S)=f(C_{x,\gamma})(S).
  \end{equation*}

Now, suppose that $S\in\CS_x$.

(i) Let $C(x)<f_C(\CS^*)=f_C(S^*)$. Taking into account \eqref{e:altma}, we have
to show that
  \begin{equation} \label{e:fifs}
f(C_{x,\gamma})(S^*)\le f(C_{x,\gamma})(S)
=\max\{f_C(S)+\gamma-C(x),F_C(S\setminus\{x\})\}.
  \end{equation}
Since $f_C(S^*)\le f_C(S)$ implies the left-hand side inequality in \eqref{e:cot}, we have
  \begin{equation} \label{e:gigs}
f_C(S^*)+\gamma-C(x)\le f(C_{x,\gamma})(S).
  \end{equation}
It follows from inequalities $C(x)<f_C(S^*)$ and
  \begin{equation*}
f_C(S^*)=\max\{C(x),F_C(S^*\setminus\{x\})\}\le f_C(S)
=\max\{C(x),F_C(S\setminus\{x\})\}
  \end{equation*}
that $F_C(S^*\setminus\{x\})\le F_C(S\setminus\{x\})\le f(C_{x,\gamma})(S)$,
which together with \eqref{e:mxy} and \eqref{e:gigs} implies \eqref{e:fifs}.

(ii) Let $C(x)=f_C(\CS^*)$ and $\CS^*=\{S^*\}$. By virtue of \eqref{e:mxy},
we have
  \begin{equation} \label{e:hops}
f(C_{x,\gamma})(S^*)=\max\{\gamma,F_C(S^*\setminus\{x\})\}.
  \end{equation}
Two cases are possible for $S\in\CS_x$ (cf.~\eqref{e:mposs}): 1) $C(x)=f_C(S)$, or
2) $f_C(S)=F_C(S\setminus\{x\})$. If case 1) holds, then $S\in\CS^*=\{S^*\}$, and so,
$S=S^*$ and inequality \eqref{e:fifs} is clear. In case 2), by virtue of \eqref{e:gorn},
we have $\gamma\le F_C(S\setminus\{x\})$, and so,
  \begin{equation} \label{e:liw}
f(C_{x,\gamma})(S)=\max\{\gamma,F_C(S\setminus\{x\})\}
=F_C(S\setminus\{x\})=f_C(S).
  \end{equation}
Thus, it follows from \eqref{e:hops}, inequality $\gamma\le C(x)$ and \eqref{e:liw} that
  \begin{equation*}
f(C_{x,\gamma})(S^*)\le\max\{C(x),F_C(S^*\setminus\{x\})\}=f_C(S^*)\le
f_C(S)=f(C_{x,\gamma})(S).
  \end{equation*}

(iii) Assume that $C(x)=f_C(\CS^*)\le F_C(S_x\setminus\{x\})$. Then, by \eqref{e:Sx},
we find $C(x)\le F_C(S\setminus\{x\})$ for $S\in\CS_x$, and so,
  \begin{equation} \label{e:ror}
f_C(S)=\max\{C(x),F_C(S\setminus\{x\})\}=F_C(S\setminus\{x\}).
  \end{equation}
The rest of this proof is as in case 2) of step~(ii).

(iv) Let $f_C(S^*)\le F_C(S_x\setminus\{x\})$. Since $x\in S^*$,
$C(x)\le f_C(S^*)\le F_C(S\setminus\{x\})$, and so, \eqref{e:ror} holds and,
by \eqref{e:49} and \eqref{e:gorn}, $f(C_{x,\gamma})(S)=f_C(S)$. Now, it
follows from \eqref{e:mxy} and inequality $\gamma\le C(x)$ that
  \begin{equation*}
f(C_{x,\gamma})(S^*)\le\max\{f_C(S^*),F_C(S^*\setminus\{x\})\}
=f_C(S^*)\le f_C(S)=f(C_{x,\gamma})(S),
  \end{equation*}
which was to be proved.
\Halmos\endproof

Now we are in a position to define the notion of the lower tolerance.

\subsection{Definition.} \label{d:loto}
Suppose $\oplus$ is a strict {\sf A}-operation on $\Rb^+$. Given $S^*\in\CS^*$
and $x\in X\setminus S^*$, the {\em lower tolerance\/} of $x$ is defined by
  \begin{equation} \label{e:ltev}
\ell_{S^*}(x)=C(x)\upsub C^-_{S^*}(x)\in\Rb^+,
  \end{equation}
where $\upsub$ is the upper subtraction for~$\oplus$.

\subsection{Theorem.} \label{t:lest} {\em
Let the\/ {\sf A}-operation $\oplus$ be strict on\/ $\Rb^+$, $S^*\in\CS^*$ and
$x\in X\setminus S^*$. Then the value $\ell_{S^*}(x)$ is well-defined,
  \begin{equation} \label{e:est}
{\pmb e}\le\ell_{S^*}(x)\le f_C(\CS^*_x)\upsub f_C(\CS^*)\equiv
\overline{\ell}_{S^*}(x),
  \end{equation}
and $\ell_{S^*}(x)\!=\!{\pmb e}$ iff\/ $C^-_{S^*}(x)\!=\!C(x)$, where
${\pmb e}$ is the neutral element with respect to~$\oplus$.

Moreover, if\/ $\oplus$ is a strict\/ {\sf A}-operation of addition on $[0,\infty)$, then
  \begin{equation} \label{e:qa}
\ell_{S^*}(x)=\left\{
  \begin{array}{ccl}
f_C(\CS^*_x)\upsub f_C(\CS^*) & \mbox{if\/} &
   F_C(S_x\setminus\{x\})\le f_C(\CS^*),\\[2pt]
C(x) & \mbox{if\/} & F_C(S_x\setminus\{x\})\ge f_C(\CS^*),
  \end{array}\right.
  \end{equation} 
and if\/ $\oplus$ is a strict\/ {\sf A}-operation of multiplication on $(0,\infty)$, then
  \begin{equation} \label{e:5.25*}
\ell_{S^*}(x)=f_C(\CS^*_x)\upsub f_C(\CS^*)=\overline{\ell}_{S^*}(x).
  \end{equation} 
}

\proof{Proof.}
That $\ell_{S^*}(x)$ is well-defined can be established along the same lines as in
the proof of Theorem~\ref{t:upto}.

Theorem~\ref{t:Cmm}(b), (in)equalities
${\pmb e}\oplus C^-_{S^*}(x)=C^-_{S^*}(x)\le C(x)$ and definition
\eqref{e:dups} imply
  \begin{equation*}
{\pmb e}\le C(x)\upsub C^-_{S^*}(x)=\ell_{S^*}(x).
  \end{equation*}

It follows from the definition of $\ell_{S^*}(x)$ and \eqref{e:ups} that
  \begin{equation} \label{e:523}
\ell_{S^*}(x)\oplus C^-_{S^*}(x)=(C(x)\upsub C^-_{S^*}(x))\oplus C^-_{S^*}(x)=C(x).
  \end{equation}
If $\ell_{S^*}(x)={\pmb e}$, then \eqref{e:523} implies $C^-_{S^*}(x)=C(x)$;
now, if $C^-_{S^*}(x)=C(x)$, then \eqref{e:523} yields equality
$\ell_{S^*}(x)\oplus C(x)={\pmb e}\oplus C(x)$, and so, taking into account the
strictness of $\oplus$ and cancelling by $C(x)$, we get $\ell_{S^*}(x)={\pmb e}$.

Let us prove the right-hand side inequality in \eqref{e:est}. We set
  \begin{equation} \label{e:uvw}
u=C(x),\quad v=f_C(\CS^*)\quad\mbox{and}\quad w=f_C(\CS^*_x).
  \end{equation}
Then it follows from \eqref{e:ltev} and the third equality in \eqref{e:Cmeq} that
  \begin{equation*}
\ell_{S^*}(x)=C(x)\upsub C^-_{S^*}(x)=u\upsub\,[(u\oplus v)\losub w].
  \end{equation*}
By virtue of Lemmas \eqref{l:luineq}(c) and \ref{l:trinv}(a), we find
  \begin{equation*}
u\upsub\,[(u\oplus v)\losub w]\le(u\oplus w)\upsub(u\oplus v)=w\upsub v,
  \end{equation*}
and it remains to note that $w\upsub v=f_C(\CS^*_x)\upsub f_C(\CS^*)%
=\overline{\ell}_{S^*}(x)$.

Let us establish \eqref{e:qa}. By virtue of \eqref{e:Cmeq},
Lemma~\ref{l:comp}(b) and \eqref{e:uvw}, we have
  \begin{equation} \label{e:bu}
C^-_{S^*}(x)=(u\oplus v)\losub w=(u\oplus v)\upsub w\quad\mbox{if}\quad
u\oplus v\ge w,
  \end{equation}
and $C^-_{S^*}(x)=0$ otherwise. Taking into account \eqref{e:Sxeq}, \eqref{e:5.4}
and the strictness of $\oplus$, we find that $u\oplus v\ge w$ iff
  \begin{equation*}
C(x)\oplus f_C(\CS^*)\ge f_C(\CS^*_x)=C(x)\oplus F_C(S_x\!\setminus\!\{x\})
\quad\mbox{iff}\quad f_C(\CS^*)\ge F_C(S_x\!\setminus\!\{x\}),
  \end{equation*}
and so, applying Lemmas \ref{l:luineq}(a) and \ref{l:trinv}(a) as well as definition
\eqref{e:ltev}, in this case we get
  \begin{equation} \label{e:ooh}
\ell_{S^*}(x)=u\upsub\,[(u\oplus v)\upsub w]=(u\oplus w)\upsub(u\oplus v)
=w\upsub v=f_C(\CS^*_x)\upsub f_C(\CS^*).
  \end{equation}

If $u\oplus v\le w$, i.e., $f_C(\CS^*)\le F_C(S_x\!\setminus\!\{x\})$, then
$\ell_{S^*}(x)=C(x)\upsub0=C(x)$, which establishes equality \eqref{e:qa}.

Finally, let us prove \eqref{e:5.25*}. If $\oplus$ is an {\sf A}-operation of
multiplication, then equalities in \eqref{e:bu} follow from \eqref{e:Cmeq} and
Lemma~\ref{l:comp}(a), (c) and hold with no restrictions on $u$, $v$ and~$w$.
This and Lemmas \ref{l:luineq}(a) and \ref{l:trinv}(a) imply equalities \eqref{e:ooh}.
\Halmos\endproof

\subsection{Examples of lower tolerances.} \label{ss:lota}
Taking into account Examples~\ref{ss:exul}(a), (c) and \ref{ss:eom}(a), (c),
\eqref{e:Cmeq}, \eqref{e:ltev}, \eqref{e:qa} and \eqref{e:bu}, given $S^*\in\CS^*$
and $x\in X\setminus S^*$, we have:

(a) if $u\oplus v=(u^p+v^p)^{1/p}$ on $\Rb^+=[0,\infty)$ with $p>0$, then
  \begin{equation*}
C^-_{S^*}(x)=\bigl(C(x)^p+f_C(\CS^*)^p-f_C(\CS^*_x)^p\bigr)^{1/p}\quad
\mbox{if}\quad F_C(S_x\!\setminus\!\{x\})\le f_C(\CS^*),
  \end{equation*} 
and $C^-_{S^*}(x)=0$ otherwise,
  \begin{equation*}
\ell_{S^*}(x)=\bigl(f_C(\CS^*_x)^p-f_C(\CS^*)^p\bigr)^{1/p}\quad\mbox{if}\quad
F_C(S_x\!\setminus\!\{x\})\le f_C(\CS^*),
  \end{equation*}
and $\ell_{S^*}(x)=C(x)$ otherwise, where (cf.~\eqref{e:Sx})
  \begin{equation*}
F_C(S_x\!\setminus\!\{x\})=\min_{S\in\CS_x}\Bigl(
\sum_{y\in S\setminus\{x\}}\!\!\!C(y)^p\Bigr)^{1/p}\quad\,\,\mbox{and}\,\,\quad
f_C(\CS^*)=\Bigl(\sum_{y\in S^*}\!C(y)^p\Bigr)^{1/p};
  \end{equation*}

(b) if $u\oplus v=u\!\cdot\! v$ on $\Rb^+=(0,\infty)$ and $C(y)>0$ for all $y\in X$,
then
  \begin{equation*}
C^-_{S^*}(x)=\frac{C(x)\!\cdot\! f_C(\CS^*)}{f_C(\CS^*_x)}\,\,\quad\mbox{and}
\,\,\quad \ell_{S^*}(x)=\frac{f_C(\CS^*_x)}{f_C(\CS^*)}=\overline{\ell}_{S^*}(x).
  \end{equation*}

\subsection{Example.} \label{ss:ses}
For a strict {\sf A}-operation of addition $\oplus$ the right-hand side inequality in
\eqref{e:est} may be strict as can be seen from the following simple example
(as well as from formula \eqref{e:qa}).

We set $X=\{x_1,x_2,x_3,x_4\}$ with $C(x_1)=C(x_2)=C(x_3)=1$ and $C(x_4)=3$,
$\CS=\{S_1,S_2\}$ with $S_1=\{x_1,x_2\}$ and $S_2=\{x_3,x_4\}$, and $\oplus=+$
(ordinary addition), and so, $f_C(S)=F_C(S)=\sum_{y\in S}C(y)$ for $S\in\CS$
(see Table~1).

\bigbreak
\begin{center}
\begin{tabular}{|c|c|c|c|c|}
\hline
$x$ & $x_1$ & $x_2$ & $x_3$ & $x_4$ \\ \hline
$C(x)$ & $1$ & $1$ & $1$ & $3$ \\ \hline\hline
$u_{S^*}(x)$ & $2$ & $2$ & & \\ \hline
$\ell_{S^*}(x)$ & & & $1$ & $2$ \\ \hline
$\overline{\ell}_{S^*}(x)$ & & & $2$ & $2$ \\ \hline
\end{tabular}\\[6pt]
Table~1
\end{center}

\bigbreak
Since $f_C(S_1)=C(x_1)+C(x_2)=2$ and $f_C(S_2)=C(x_3)+C(x_4)=4$, we find
$\CS^*=\{S^*\}$ with $S^*=S_1$ and $f_C(\CS^*)=f_C(S^*)=2$. Noting that
$\CS_{-x}=\{S_2\}$ and $f_C(\CS^*_{-x})=f_C(S_2)=4$ for $x=x_1,x_2$ and
taking into account Examples \ref{ss:exul}(a) and \ref{ss:ust}(a) with $p=1$, we have:
  \begin{equation*}
u_{S^*}(x)=f_C(\CS^*_{-x})-f_C(\CS^*)=4-2=2\quad\,\,\mbox{if}
\,\,\quad x=x_1,\,x_2.
  \end{equation*}
Now, since $S_{x_3}=S_{x_4}=S_2$, $S_{x_3}\setminus\{x_3\}=\{x_4\}$ and
$S_{x_4}\setminus\{x_4\}=\{x_3\}$, and so,
  \begin{equation*}
F_C(S_{x_3}\!\setminus\!\{x_3\})\!=\!C(x_4)\!=\!3\!>\!f_C(\CS^*)
\quad\mbox{and}\quad
F_C(S_{x_4}\!\setminus\!\{x_4\})\!=\!C(x_3)\!=\!1\!<\!f_C(\CS^*),
  \end{equation*}
it follows from \eqref{e:qa} and Example~\ref{ss:ses}(a) that
  \begin{equation*}
\ell_{S^*}(x_3)=C(x_3)=1\quad\mbox{and}\quad \ell_{S^*}(x_4)%
=f_C(\CS^*_{x_4})-f_C(\CS^*)=4-2=2.
  \end{equation*}
On the other hand, since $\CS_x=\{S_2\}$ for $x=x_3,x_4$, we find
(cf.~\eqref{e:est})
  \begin{equation*}
\overline{\ell}_{S^*}(x)=f_C(\CS^*_x)-f_C(\CS^*)=4-2=2\quad\mbox{if}
\quad x=x_3,x_4.
  \end{equation*}

\subsection{Remark.} \label{r:costs}
The quantity $\overline{\ell}_{S^*}(x)$ from \eqref{e:est}, which may be called the
{\em extended lower tolerance\/} of~$x$, differs from $\ell_{S^*}(x)$ in the following
way. Because $C^-_{S^*}(x)\in\Rb^+$ and $C^-_{S^*}(x)\le C(x)$, elements of
the interval $[C^-_{S^*}(x),C(x)]$ are still {\em costs\/} (i.e., nonnegative).
Moreover, by virtue of \eqref{e:ltev} and \eqref{e:ups},
$\ell_{S^*}(x)\oplus C^-_{S^*}(x)=C(x)$, and so,
$C^-_{S^*}(x)=C(x)\upsub\ell_{S^*}(x)$ is restored from $\ell_{S^*}(x)$ as the
lowest possible stability cost of~$x$. From this point of view the quantity
$\widetilde{C}^-_{S^*}(x)=C(x)\upsub\,\overline{\ell}_{S^*}(x)$ may not be defined
as a cost. For instance, in Example~\ref{ss:ses} we have
$C^-_{S^*}(x_3)=1-1=0$, whereas $\widetilde{C}^-_{S^*}(x_3)=1-2=-1\notin\Rb^+$.

\section{Tolerance functions.} \label{s:TF}

Throughout this section we assume that the \DOP$(X,\CS,\oplus,C)$ of the form
\eqref{e:DO} is given, $\oplus$ is {\em strict\/} and $\CS^*=\CS^*_C$ is
from~\eqref{e:S*}.

We are going to define a function $T_C$ on $X$, which is an {\em invariant\/} of the
\DOP{} under consideration in the sense that it is {\em independent\/} of optimal
solutions~$S^*\in\CS^*$.

If $\oplus$ is an {\sf A}-operation of addition on $\Rb^+=[0,\infty)$ and $u\ge0$, then
we set $u^{-1}=-u$, and if $\oplus$ is an {\sf A}-operation of multiplication on
$\Rb^+=(0,\infty)$ with the neutral element ${\pmb e}\in\Rb^+$ and $u>0$, then,
taking into account Lemma~\ref{l:comp}(c), we set $u^{-1}={\pmb e}\upsub u$.
It is to be noted that in the latter case we have
  \begin{equation} \label{e:-1}
u^{-1}\oplus u={\pmb e},\quad{\pmb e}^{-1}={\pmb e},\quad
(u^{-1})^{-1}=u,\quad\mbox{and}\quad\mbox{$u>{\pmb e}$ \,\,iff \,
$u^{-1}<{\pmb e}$.}
  \end{equation}
The last three properties in \eqref{e:-1} also hold in the former case with ${\pmb e}=0$.

\subsection{Definition.} \label{d:tofu}
Given $S^*\in\CS^*$ and $x\in X$, we set
  \begin{equation} \label{e:tofu}
T_C(x)=\left\{
  \begin{array}{ccl}
u_{S^*}(x) & \mbox{if} & x\in S^*,\\[2pt]
(\ell_{S^*}(x))^{-1} & \mbox{if} & x\in X\setminus S^*.
  \end{array}\right.
  \end{equation}
The function $T_C$ on $X$ is said to be the {\em tolerance function\/} of the
\DOP$(X,\CS,\oplus,C)$. Replacing $\ell_{S^*}(x)$ by $\overline{\ell}_{S^*}(x)$ in
the second line of \eqref{e:tofu} we get the notion of the {\em extended tolerance
function}, denoted by $\overline T_C(x)$. Both tolerance functions assume their values
in $\Rb$ if $\oplus$ is an {\sf A}-operation of addition and in $\Rb^+=(0,\infty)$
if $\oplus$ is an {\sf A}-operation of multiplication. Note also that $T_C(x)$ can be
represented as
  \begin{equation*}
T_C(x)=u_{S^*}(x)\chi_{S^*}(x)+(\ell_{S^*}(x))^{-1}\chi_{X\setminus S^*}(x),
\qquad x\in X,
  \end{equation*}
where, given $Y\subset X$, $\chi_Y$ is the characteristic function of the set $Y$
(i.e., $\chi_Y(x)=1$ if $x\in Y$ and $\chi_Y(x)=0$ if $x\in X\setminus Y$).

\smallbreak
The correctness of this definition is justified by the following

\subsection{Theorem.} \label{t:corTF} {\em
The tolerance function $T_C$ on $X$ is well-defined, independent of optimal solutions
$S^*\in\CS^*$ and has the following properties\/{\rm:}
  \begin{equation} \label{e:Ua}
T_C|_{S^*}(\cdot)\!=\!u_{S^*}(\cdot)\!\ge\!{\pmb e}\quad\mbox{and}\quad
T_C|_{X\setminus S^*}(\cdot)\!=\!(\ell_{S^*}(\cdot))^{-1}\!\le\!{\pmb e}\quad
\mbox{for all \,$S^*\in\CS^*$,}
  \end{equation}
where $T_C|_Y(\cdot)$ denotes the restriction of\/ $T_C$ to the set $Y\subset X$.
}

\medbreak
In order to prove Theorem~\ref{t:corTF}, we need a lemma, which is of interest
in its own.

\subsection{Lemma.} \label{l:KK} {\em
Given $S_1^*,S_2^*\in\CS^*$ and $x\in X$, we have\/{\rm:}
\par\smallbreak
{\rm(a)} if\/ $x\in S_1^*\cap S_2^*$, then $u_{S_1^*}(x)=u_{S_2^*}(x)$
{\rm(}and $\ge{\pmb e}${\rm);}
\par\smallbreak
{\rm(b)} if\/ $x\in(X\setminus S_1^*)\cap(X\setminus S_2^*)$, then
$\ell_{S_1^*}(x)=\ell_{S_2^*}(x)$ {\rm(}and $\ge{\pmb e}${\rm);}
\par\smallbreak
{\rm(c)} if\/ $x\in S_1^*\setminus S_2^*$, then
$u_{S_1^*}(x)={\pmb e}=\ell_{S_2^*}(x)$.
}
\medbreak
\proof{Proof.}
(a) By virtue of \eqref{e:ucal}, we have
  \begin{equation*}
u_{S_1^*}(x)=f_C(\CS^*_{-x})\upsub f_C(\CS^*)=u_{S_2^*}(x).
  \end{equation*}

(b) It follows from \eqref{e:ltev} and \eqref{e:Cmeq} that
  \begin{align*}
\ell_{S_1^*}(x)&=C(x)\upsub C^-_{S_1^*}(x)=C(x)\upsub\bigl(
  [C(x)\oplus f_C(\CS^*)]\losub f_C(\CS^*_x)\bigr)=\\[2pt]
&=C(x)\upsub C^-_{S_2^*}(x)=\ell_{S_2^*}(x).
  \end{align*}

(c) Since $x\in S_1^*$ iff $S_1^*\in\CS_x$, and $x\in X\setminus S_2^*$ iff
$S_2^*\in\CS_{-x}$, \eqref{e:CS} and \eqref{e:ov} imply
$f_C(\CS^*_{-x})\le f_C(S_2^*)=f_C(\CS^*)$. Now, it follows from
Theorem \ref{t:upto} and Lemma~\ref{l:fivar}(a) that
  \begin{equation*}
{\pmb e}\le u_{S_1^*}(x)=f_C(\CS^*_{-x})\upsub f_C(\CS^*)\le
f_C(\CS^*)\upsub f_C(\CS^*)={\pmb e},
  \end{equation*}
and so, $u_{S_1^*}(x)={\pmb e}$. At the same time, by virtue of \eqref{e:CS},
we find $f_C(\CS^*_x)\le f_C(S_1^*)=f_C(\CS^*)$, and so, Lemma~\ref{l:scnde}(b),
the strictness of $\oplus$ and \eqref{e:wvv} yield
  \begin{equation*}
[C(x)\oplus f_C(\CS^*)]\losub f_C(\CS^*_x)\ge
[C(x)\oplus f_C(\CS^*)]\losub f_C(\CS^*)=C(x).
  \end{equation*}
Now, it follows from Theorem~\ref{t:lest}, \eqref{e:ltev}, \eqref{e:Cmeq} and
Lemma~\ref{l:scnde}(a) that
  \begin{equation*}
{\pmb e}\le\ell_{S_2^*}(x)=C(x)\upsub\bigl([C(x)\oplus f_C(\CS^*)]
\losub f_C(\CS^*_x)\bigr)\le C(x)\upsub C(x)={\pmb e},
  \end{equation*}
which implies the equality $\ell_{S_2^*}(x)={\pmb e}$.
\Halmos\endproof

Now we are in a position to prove Theorem~\ref{t:corTF}.
\medbreak
\proof{Proof of Theorem\/~\ref{t:corTF}.}
Formally, the tolerance function, as it is defined in \eqref{e:tofu}, depends on
$S^*\in\CS^*$, and so, we temporarily write $T_{C,S^*}$ in place of~$T_C$.
Only the case $|\CS^*|\ge2$ is to be considered, and so, assuming that
$S_1^*,S_2^*\in\CS^*$ are arbitrarily chosen, let us show that
$T_{C,S_1^*}(x)=T_{C,S_2^*}(x)$ for all $x\in X$.

In fact, we have the following decomposition of the set $X$:
  \begin{equation*}
X=[S_1^*\cap S_2^*]\cup[(X\setminus S_1^*)\cap(X\setminus S_2^*)]\cup
[(S_1^*\cup S_2^*)\setminus(S_1^*\cap S_2^*)],
  \end{equation*}
where, by de\,Morgan's laws, $(X\!\setminus\! S_1^*)\cap(X\!\setminus\! S_2^*)=%
X\!\setminus\!(S_1^*\cup S_2^*)$, and the sets on the right in square brackets
are pairwise disjoint. Suppose $x\in X$. If $x\in S_1^*\cap S_2^*$, then \eqref{e:tofu}
and Lemma~\ref{l:KK}(a) imply
  \begin{equation*}
T_{C,S_1^*}(x)=u_{S_1^*}(x)=u_{S_2^*}(x)=T_{C,S_2^*}(x).
  \end{equation*}
If $x\in(X\setminus S_1^*)\cap(X\setminus S_2^*)$, then it follows from \eqref{e:tofu}
and Lemma~\ref{l:KK}(b) that
  \begin{equation*}
T_{C,S_1^*}(x)=(\ell_{S_1^*}(x))^{-1}=(\ell_{S_2^*}(x))^{-1}=T_{C,S_2^*}(x).
  \end{equation*}
Now, assume that
  $
x\in(S_1^*\cup S_2^*)\setminus(S_1^*\cap S_2^*)
=(S_1^*\setminus S_2^*)\cup(S_2^*\setminus S_1^*)
  $
and, to be more specific, $x\in S_1^*\setminus S_2^*$. Then, by virtue of
\eqref{e:tofu}, Lemma~\ref{l:KK}(c) and \eqref{e:-1}, we find
  \begin{equation*}
T_{C,S_1^*}(x)=u_{S_1^*}(x)={\pmb e}={\pmb e}^{-1}=
(\ell_{S_2^*}(x))^{-1}=T_{C,S_2^*}(x).
  \end{equation*}

Now it is correct to set back $T_C(x)=T_{C,S^*}(x)$ for all $x\in X$ and $S^*\in\CS^*$,
which together with \eqref{e:tofu} implies
$T_C|_{S^*}(x)=T_C(x)=u_{S^*}(x)$ if $x\in S^*$ and
$T_C|_{X\setminus S^*}(x)=T_C(x)=(\ell_{S^*}(x))^{-1}$ if $x\in X\setminus S^*$.
This completes the proof.
\Halmos\endproof

\subsection{Remark.} \label{r:jr}
An assertion similar to Theorem~\ref{t:corTF} holds for the extended tolerance
function $\overline T_C$ with obvious modifications. In fact, it is to be noted only that
in Lemma~\ref{l:KK}(b) we have, by virtue of \eqref{e:est},
$\overline{\ell}_{S^*_1}(x)=f_C(\CS^*_x)\upsub f_C(\CS^*)=\overline{\ell}_{S_2^*}(x)$,
and under conditions of Lemma~\ref{l:KK}(c), we get
  \begin{equation*}
{\pmb e}\le\overline{\ell}_{S_2^*}(x)=f_C(\CS^*_x)\upsub f_C(\CS^*)\le
f_C(\CS^*)\upsub f_C(\CS^*)={\pmb e}.
  \end{equation*}

\subsection{Example.} \label{ss:just}
Let us illustrate Theorem~\ref{t:corTF} by the following example of small cardinality
$|X|$ and the simplest possible {\sf A}-operation $\oplus=+$ (so that calculations
are not cumbersome and all the details can be clearly seen).

Let $X=\{x_1,x_2,x_3,x_4,x_5,x_6\}$ with $C(x_1)=C(x_2)=C(x_3)=2$, $C(x_4)=1$,
$C(x_5)=3$ and \mbox{$C(x_6)=5$}, and $\CS=\{S_1,S_2,S_3\}$ with
$S_1=\{x_1,x_2,x_3\}$,
$S_2=\{x_2,x_4,x_5\}$ and $S_3=\{x_1,x_4,x_6\}$. Since the objective function
is of the form $f_C(S)=\sum_{y\in S}C(y)$, $S\in\CS$, we find
$f_C(S_1)=f_C(S_2)=6$ and $f_C(S_3)=8$, and so, $\CS^*=\{S_1^*,S_2^*\}$
with $S_1^*=S_1$ and $S_2^*=S_2$. The corresponding values of upper and lower
tolerances and the tolerance function are presented in the following Table~2:

\bigbreak
\begin{center}
\begin{tabular}{|c|c|c|c|c|c|c|}
\hline
$x$ & $x_1$ & $x_2$ & $x_3$ & $x_4$ & $x_5$ & $x_6$\\ \hline
$C(x)$ & $2$ & $2$ & $2$ & $1$ & $3$ & $5$\\ \hline
$S_1$ & $*$ & $*$ & $*$ & & & \\ \hline
$S_2$ & & $*$ & & $*$ & $*$ & \\ \hline
$S_3$ & $*$ & & & $*$ & & $*$\\ \hline\hline
$u_{S^*_1}(x)$ & $0$ & $2$ & $0$ & & &\\ \hline
$\ell_{S^*_1}(x)$ & & & & $0$ & $0$ & $2$\\ \hline\hline
$T_C(x)$ & $0$ & $2$ & $0$ & $0$ & $0$ & $-2$ \\ \hline\hline
$u_{S^*_2}(x)$ & & $2$ & & $0$ & $0$ &\\ \hline
$\ell_{S^*_2}(x)$ & $0$ & & $0$ & & & $2$\\ \hline
\end{tabular}\\[6pt]
Table~2
\end{center}

\bigbreak
Given $i\in\{1,2,3\}$ and $j\in\{1,2,3,4,5,6\}$, we put $*$ in row $S_i$ and column $x_j$
provided $x_j\in S_i$. Then setting $S^*=S_1^*$ we calculate the values
$u_{S_1^*}(x)$ for $x\in S_1^*$ and $\ell_{S_1^*}(x)$ for $x\in X\setminus S_1^*$
in accordance with Theorems~\ref{t:upto} and \ref{t:lest} and Lemma~\ref{l:KK}(c).
By virtue of \eqref{e:tofu}, we form the tolerance function
  \begin{equation*}
T_C(x)=\bigl(T_C(x_1),T_C(x_2),T_C(x_3),T_C(x_4),T_C(x_5),T_C(x_6)\bigr)
=(0,2,0,0,0,-2).
  \end{equation*}
Now, making use of Theorem~\ref{t:corTF} and taking into account row $S_2$, 
in which elements of the optimal solution $S_2^*$ and outside of it are marked, we
extract from the vector $T_C(x)$ the corresponding values of $u_{S_2^*}(x)$
for $x\in S_2^*$ and $\ell_{S_2^*}(x)$ for $x\in X\setminus S_2^*$.

Thus, Theorem~\ref{t:corTF} says that only upper and lower tolerances with respect
to any fixed optimal solution to the \DOP{} under consideration are to be calculated,
the other tolerances being determined uniquely via the tolerance function.

\subsection{Convention.} \label{ss:gez}
In what follows we assume that $C(x)>0$ for all $x\in X$.

\smallbreak
Definition \ref{d:tofu} is also motivated by the fact that the set of optimal solutions
$\CS^*$ to the discrete optimization problem \eqref{e:DO} can be characterized
by means of the tolerance function(s) in the following way.

\subsection{Theorem.} \label{t:TCxe} {\em
{\rm(a)} $\{x\in X:T_C(x)={\pmb e}\}=(\cup\CS^*)\setminus(\cap\CS^*)$.
  \begin{itemize}
\item[{\rm(b)}] $\{x\in X:T_C(x)>{\pmb e}\}=\cap\CS^*$.
\item[{\rm(c)}] $\{x\in X:T_C(x)\ge{\pmb e}\}=\cup\CS^*$.
\item[{\rm(d)}] $\{x\in X:T_C(x)<{\pmb e}\}=X\setminus(\cup\CS^*)$.
\item[{\rm(e)}] $\{x\in X:T_C(x)\le{\pmb e}\}=X\setminus(\cap\CS^*)$.
  \end{itemize}
}
\medbreak
\proof{Proof.}
(a)\,($\supset$) If $x\in(\cup\CS^*)\setminus(\cap\CS^*)$, then there exist two
optimal trajectories $S_1^*,S_2^*\in\CS^*$ such that $x\in S_1^*\setminus S_2^*$,
and so, by \eqref{e:tofu} and Lemma~\ref{l:KK}(c), we get
  \begin{equation*}
T_C(x)=u_{S_1^*}(x)=(\ell_{S_2^*}(x))^{-1}={\pmb e}.
  \end{equation*}

(a)\,($\subset$) Suppose $x\in X$ and $T_C(x)={\pmb e}$. Let us fix an $S^*\in\CS^*$.
If $x\in S^*$, then, by virtue of \eqref{e:tofu} and \eqref{e:ucal}, we have
  \begin{equation*}
{\pmb e}=T_C(x)=u_{S^*}(x)=f_C(\CS^*_{-x})\upsub f_C(\CS^*),
  \end{equation*}
which implies $f_C(\CS^*_{-x})=f_C(\CS^*)$ (cf.~\eqref{e:ups}). Taking into account
\eqref{e:CS}, we may choose $S_1\in\CS_{-x}$ such that $f_C(S_1)=f_C(\CS^*)$,
whence $S_1\in\CS^*$. Since $x\notin S_1$ and $x\in S^*$, we get
$x\in S^*\setminus S_1\subset(\cup\CS^*)\setminus(\cap\CS^*)$.

Now, assume that $x\notin S^*$. By virtue of \eqref{e:tofu}, \eqref{e:qa}
and \eqref{e:5.25*}, we claim that
  \begin{equation} \label{e:A*}
{\pmb e}={\pmb e}^{-1}=(T_C(x))^{-1}=\ell_{S^*}(x)
=f_C(\CS^*_x)\upsub f_C(\CS^*).
  \end{equation}
Only the last equality in \eqref{e:A*} is to be justified in the case when $\oplus$
is an {\sf A}-operation of addition on $[0,\infty)$: on the contrary, if this is not so,
then \eqref{e:qa} and \eqref{e:tofu} imply
  \begin{equation*}
C(x)=\ell_{S^*}(x)=(T_C(x))^{-1}={\pmb e}^{-1}={\pmb e}=0,
  \end{equation*}
which contradicts our convention~\ref{ss:gez}. Making use of \eqref{e:A*}, we get
$f_C(\CS^*_x)=f_C(\CS^*)$, and so, by virtue of \eqref{e:CS}, there exists
$S_2\in\CS_x$ such that $f_C(S_2)=f_C(\CS^*)$ and, hence, $S_2\in\CS^*$.
Since $x\in S_2$ and $x\notin S^*$, we find
$x\in S_2\setminus S^*\subset(\cup\CS^*)\setminus(\cap\CS^*)$.

\smallbreak
(b)\,($\subset$) Let $x\in X$ and $T_C(x)>{\pmb e}$. We claim that
$x\in\cap\CS^*$. On the contrary, assume that $x\notin\cap\CS^*$, and so,
$x\notin S^*$ for some $S^*\in\CS^*$. By Theorem~\ref{t:lest},  we have
$\ell_{S^*}(x)\ge{\pmb e}$, and so, \eqref{e:tofu} gives
$T_C(x)=(\ell_{S^*}(x))^{-1}\le{\pmb e}$, which is a contradiction.

\smallbreak
(b)\,($\supset$) Let $x\in\cap\CS^*$. Choose an $S^*\in\CS^*$. Since $x\in S^*$,
\eqref{e:tofu} and \eqref{e:ucal} imply $T_C(x)=u_{S^*}(x)\ge{\pmb e}$. Taking
into account item (a), we infer that $T_C(x)>{\pmb e}$.

\smallbreak
(c) is a consequence of (a) and (b): $T_C(x)\ge{\pmb e}$ iff $T_C(x)={\pmb e}$ or
$T_C(x)>{\pmb e}$, i.e., iff $x\in(\cup\CS^*)\setminus(\cap\CS^*)$ or
$x\in\cap\CS^*$.

\smallbreak
(d) follows immediately from (c):
  \begin{equation*}
\{x\in X:T_C(x)<{\pmb e}\}=X\setminus\{x\in X:T_C(x)\ge{\pmb e}\}
=X\setminus(\cup\CS^*).
  \end{equation*}

(e) is a straightforward consequence of (b).
\Halmos\endproof

Tolerance functions can be effectively applied for the characterization of uniqueness and
nonuniqueness of optimal trajectories:

\subsection{Corollary.} \label{c:uni} {\em
{\rm(a)} The\/ {\rm \DOP{} \eqref{e:DO}} admits a unique optimal solution\/
{\rm(}i.e., \mbox{$|\CS^*|=1$}{\rm)} if and only if\/ $T_C(\cdot)\ne{\pmb e}$ on $X$.

{\rm(b)} $|\CS^*|\ge2$ iff\/ $T_C(x)={\pmb e}$ for some $x\in X$.
}
\medbreak
\proof{Proof.}
(a) By virtue of Theorem~\ref{t:TCxe}(a), $|\CS^*|=1$ iff $\cup\CS^*=\cap\CS^*$
iff $(\cup\CS^*)\setminus(\cap\CS^*)=\varnothing$ iff $T_C(x)\ne{\pmb e}$
for all $x\in X$.

\smallbreak
(b) is simply the negation of item~(a).
\Halmos\endproof

At the end of this section we are going to establish certain relationships between
the values of $T_C$ on $S^*\in\CS^*$ and on $X\setminus S^*$.

\subsection{Covering trajectories.} \label{ss:ct}
Given $Y\subset X$, it is convenient to introduce the collection $\CS_c(Y)$
(possibly, empty) of those trajectories $S\in\CS$, which cover the set~$Y$:
  \begin{equation*}
\CS_c(Y)=\{S\in\CS:\mbox{$Y\subset S$ and $S\ne Y$}\}.
  \end{equation*}
It is to be noted that $\CS_c(Y)=\varnothing$ iff $Y\setminus S\ne\varnothing$
for all $S\in\CS$ with $S\ne Y$, and if $\CS_c(Y)\ne\varnothing$, then
$S\notin\CS_c(Y)$ iff $Y\setminus S\ne\varnothing$.

We say that the set of trajectories $\CS$ {\em consists of nonembedded sets\/}
provided that $\CS_c(S)=\varnothing$ for all $S\in\CS$. In other words
(cf.\ \cite[Theorem~1]{JOGO12}), the last condition is equivalent to saying that
$S_1\setminus S_2\ne\varnothing$ for all $S_1,S_2\in\CS$, $S_1\ne S_2$.
For instance, in Examples \ref{ss:ses} and \ref{ss:just} collections of trajectories
$\CS$ consist of nonembedded sets.

\subsection{Theorem.} \label{t:minv} {\em
Assume that $C(x)>0$ for all $x\in X$ if $\oplus$ is an\/ {\sf A}-operation of addition
on $\Rb^+=[0,\infty)$ and $C(x)\ge{\pmb e}$ for all $x\in X$ if $\oplus$ is an\/
{\sf A}-operation of multiplication on $\Rb^+=(0,\infty)$. If $S^*\in\CS^*$ is the
unique optimal solution to the\/ {\rm \DOP~\eqref{e:DO}}, then we have the inequalities
  \begin{equation} \label{e:feq}
\min_{y\in X\setminus S^*}(T_C(y))^{-1}\le
\min_{y\in X\setminus S^*}(\overline T_C(y))^{-1}\le
\min_{x\in S^*}T_C(x)
  \end{equation}
and
  \begin{equation} \label{e:zeq}
\min_{x\in S^*}T_C(x)\le\min_{y\in X\setminus[S^*\cup(\cup\CS_c(S^*))]}
(\overline T_C(y))^{-1}.
  \end{equation}
}

\proof{Proof.}
1. We begin by proving the right-hand side inequality in \eqref{e:feq} (the left-hand
side inequality in \eqref{e:feq} is always valid by virtue of \eqref{e:tofu} and
\eqref{e:est}). Given $x\in S^*$, it follows from \eqref{e:tofu} and \eqref{e:ucal} that
  \begin{equation*}
T_C(x)=u_{S^*}(x)=f_C(\CS^*_{-x})\upsub f_C(\CS^*),
  \end{equation*}
and so, by virtue of \eqref{e:CS}, there exists $S_1\in\CS_{-x}$ such that
  \begin{equation*}
T_C(x)=f_C(S_1)\upsub f_C(\CS^*).
  \end{equation*}
Also, it follows from Theorem~\ref{t:TCxe}(c) and Corollary~\ref{c:uni}(a) that
$T_C(x)>{\pmb e}$, and so, by \eqref{e:ups} and (A.3$_{\mbox{\footnotesize\rm s}}$),
  \begin{equation} \label{e:S1}
f_C(S_1)=T_C(x)\oplus f_C(\CS^*)>{\pmb e}\oplus f_C(\CS^*)=f_C(\CS^*),
  \end{equation}
i.e., $S_1\notin\CS^*$. We claim that $S_1\setminus S^*\ne\varnothing$. On the
contrary, assume that $S_1\setminus S^*=\varnothing$, and so, $S_1\subset S^*$,
say, $S_1=\{x_1,\dots,x_n\}$ and $S^*=\{x_1,\dots,x_n,y_1,\dots,y_m\}$. Then
  \begin{align*}
f_C(S_1)&=\bigoplus_{i=1}^nC(x_i)=\biggl(\bigoplus_{i=1}^nC(x_i)\biggr)\oplus
  \biggl(\bigoplus_{j=1}^m{\pmb e}\biggr)\le\\
&\le\biggl(\bigoplus_{i=1}^nC(x_i)\biggr)\oplus
  \biggl(\bigoplus_{j=1}^mC(y_j)\biggr)=f_C(S^*)=f_C(\CS^*),
  \end{align*}
which contradicts to inequality \eqref{e:S1}. Now, pick $y_0\in S_1\setminus S^*$.
Then $y_0\in S_1$ and $y_0\notin S^*$ or, in other words, $S_1\in\CS_{y_0}$ and
$y_0\in X\setminus S^*$. By virtue of \eqref{e:tofu}, \eqref{e:est}, \eqref{e:CS}
and Lemma~\ref{l:fivar}(a), we get
  \begin{align*}
\min_{y\in X\setminus S^*}(\overline T_C(y))^{-1}&\le
  (\overline T_C(y_0))^{-1}=\overline{\ell}_{S^*}(y_0)=
  f_C(\CS^*_{y_0})\upsub f_C(\CS^*)\le\\
&\le f_C(S_1)\upsub f_C(\CS^*)=T_C(x),
  \end{align*}
from which the right-hand side inequality in \eqref{e:feq} follows if we take into account
the arbitrariness of $x\in S^*$.

2. Now we establish inequality \eqref{e:zeq}. Let
$y\in X\setminus[S^*\cup(\cup\CS_c(S^*))]$. Since $y\in X\setminus S^*$, it
follows from \eqref{e:tofu} and \eqref{e:est} that
  \begin{equation*}
(\overline T_C(y))^{-1}=\overline\ell_{S^*}(y)=f_C(\CS^*_y)\upsub f_C(\CS^*),
  \end{equation*}
and so, by \eqref{e:CS}, there exists $S_2\in\CS_y$ (i.e., $S_2\in\CS$ and $y\in S_2$)
such that
  \begin{equation*}
(\overline T_C(y))^{-1}=f_C(S_2)\upsub f_C(\CS^*).
  \end{equation*}
By Theorem~\ref{t:TCxe}(e) and Corollary~\ref{c:uni}(a), $T_C(y)<{\pmb e}$,
and so, \eqref{e:-1}, \eqref{e:ups} and (A.3$_{\mbox{\footnotesize\rm s}}$) yield
that $(\overline T_C(y))^{-1}>{\pmb e}$ and
  \begin{equation*}
f_C(S_2)=(\overline T_C(y))^{-1}\oplus f_C(\CS^*)>{\pmb e}\oplus f_C(\CS^*)
=f_C(\CS^*)
  \end{equation*}
i.e., $S_2\notin\CS^*$ implying $S_2\ne S^*$.
We claim that $S^*\setminus S_2\ne\varnothing$. There are two possibilities: either
$\CS_c(S^*)=\varnothing$ or $\CS_c(S^*)\ne\varnothing$. If $\CS_c(S^*)=\varnothing$,
then no $S\in\CS$, $S\ne S^*$, covers $S^*$, and so,
$S^*\setminus S_2\ne\varnothing$. Assume that $\CS_c(S^*)\ne\varnothing$.
Because $y\notin\cup\CS_c(S^*)$, we have $y\notin S$ for all $S\in\CS_c(S^*)$
(i.e., for all $S\in\CS$ such that $S^*\subset S$ and $S\ne S^*$). Taking into account
that $y\in S_2$ and $S_2\ne S^*$, we find $S_2\notin\CS_c(S^*)$, and so,
$S_2$ does not cover $S^*$ and $S^*\setminus S_2\ne\varnothing$.
Now, choose an $x_0\in S^*\setminus S_2$. This gives $x_0\in S^*$ and
$S_2\in\CS_{-x_0}$, and so, applying \eqref{e:tofu}, \eqref{e:ucal}, \eqref{e:CS}
and Lemma~\ref{l:fivar}(a), we find
  \begin{align*}
\min_{x\in S^*}T_C(x)&\le T_C(x_0)=u_{S^*}(x_0)=
  f_C(\CS^*_{-x_0})\upsub f_C(\CS^*)\le\\[2pt]
&\le f_C(S_2)\upsub f_C(\CS^*)=(\overline T_C(y))^{-1},
  \end{align*}
and it remains to take into account the arbitrariness of $y$ as above.
\Halmos\endproof 

\subsection{Corollary.} \label{c:mv} {\em
Under the assumptions of Theorem\/~{\rm\ref{t:minv}}, if the set of trajectories
$\CS$ consists of nonembedded sets, then
  \begin{equation*}
\min_{y\in X\setminus S^*}(T_C(y))^{-1}\le
\min_{y\in X\setminus S^*}(\overline T_C(y))^{-1}=\min_{x\in S^*}T_C(x),
  \end{equation*}
and in the case of the\/ {\sf A}-operation of multiplication $\oplus$ on $\Rb^+=(0,\infty)$
the inequality $\le$ above turns out to be the equality. 
}
\medbreak
\proof{Proof.}
Since $\CS_c(S^*)=\varnothing$, the (in)equalities follow from \eqref{e:feq}
and \eqref{e:zeq}.
\Halmos\endproof

A discussion of issues as in Theorem~\ref{t:minv} is presented in
\cite[Sections 5--7]{JOGO12} in the case $\oplus=+$ (cf.\ also \cite{GS}).
Also, Examples 1 and 2 from \cite{JOGO12} show that inequalities \eqref{e:feq}
and \eqref{e:zeq} may be strict.

\bigbreak
\section*{Acknowledgments.}
The authors are supported by LATNA Laboratory, NRU HSE, RF government grant,
ag.~11.G34.31.0057. The authors are grateful to Boris Goldengorin for several
insightful discussions on the results of this paper.




\bibliographystyle{ormsv080}

\end{document}